\documentclass{amsart}
\usepackage[USenglish]{babel}
\usepackage{color}
\usepackage{graphicx}

\usepackage{amssymb}
\usepackage{amsmath}

\newcommand{\R}{\mathrm{I\!R\!}}

\newtheorem{theorem}{Theorem}[section]
\newtheorem{lemma}[theorem]{Lemma}
\newtheorem{proposition}[theorem]{Proposition}

\newtheorem{remark}[theorem]{Remark}

\numberwithin{equation}{section}

\title[A model for tumor growth and chemotherapy]{An analysis of a mathematical model describing the growth of a tumor treated with chemotherapy}

\author{Anderson L. A. de Araujo}
\address{Departamento de Matem\'atica, Universidade Federal de Vi\c{c}osa, Vi\c{c}osa, MG, Brazil}
\email{anderson.araujo@ufv.br}

\author{Artur C. Fassoni}
\address{Instituto de Matem\'atica e Computa\c{c}\~ao, Universidade Federal de Itajub\'a, Itajub\'a, MG, Brazil}
\email{fassoni@unifei.edu.br}

\author{Lu\'is F. Salvino }
\address{Instituto de Ci\^encias Exatas, Universidade Federal de Vi\c{c}osa - Campus Florestal, Florestal, MG, Brazil}
\email{luisfsalvino@gmail.com}



\begin{document}

\begin{abstract}
We present a mathematical analysis of a mixed ODE-PDE model describing the spatial distribution and temporal evolution of tumor and normal cells within a tissue subject to the effects of a chemotherapeutic drug. The model assumes that the influx of chemotherapy is restricted to a limited region of the tissue, mimicking a blood vessel passing transversely. We provide results on the existence and uniqueness of the model solution and numerical simulations illustrating different model behaviors.
\end{abstract}

\maketitle


\par
\section{Introduction}

The main objective of this work is to perform a rigorous mathematical analysis of a system of nonlinear partial differential equations corresponding to a generalization of a mathematical model describing the growth of a tumor proposed in \cite{Fassoni}.

To describe the model,
let $\Omega \subset \R^{2}$, be an open and bounded set; let also  $0< T< \infty$ be a given final time  of interest and denote $t$ the times between $[0,T]$ and
$Q=\Omega\times(0,T)$, the space-time cylinder and
$\bar{\Gamma}=\partial \Omega\times(0,T)$,  the space-time boundary. Then, the system of equations we are considering is the following:\begin{equation}
\label{0riginalEquations}
\left\{
\begin{array}{lcl}
\displaystyle	
\frac{\partial N}{\partial t} = r_N - \mu_N N - \beta_1 N A - \alpha_N\gamma_N D N, & \textup{in}& Q,
\vspace{0.2cm}
\\
\displaystyle	
\frac{\partial A}{\partial t} = r_A A\left(1-\frac{A}{k_A}\right)-(\mu_A+\epsilon_A)A - \alpha_A\gamma_A D A, &\textup{in}& Q,
\vspace{0.2cm}
\\
\displaystyle	
\frac{\partial D}{\partial t} = \sigma \Delta D + \mu\chi_{\omega} - \gamma_A D A - \gamma_N D N - \tau D, &\textup{in}& Q,
\vspace{0.2cm}
\\
\displaystyle \frac{\partial D}{\partial \eta} =0, &\textup{on}& \Gamma,
\vspace{0.2cm}
\\
\displaystyle N(\cdot,0) = N_{0} (\cdot), A(\cdot,0) = A_{0}(\cdot), D(\cdot,0) = D_{0} (\cdot), &\textup{in}& \Omega.
\end{array}
\right.
\end{equation}

In \cite{Fassoni}, Fassoni studied an ODE system corresponding to system \eqref{0riginalEquations} in a spatially homogeneous setting. Such model describes the growth of a tumor and its effect on the normal tissue, the tissue response to the tumor and the application of chemotherapeutic treatments, without spatial heterogeneity. The aim of the authors was to understand the phenomena of cancer onset and treatment as transitions between different basins of attraction of the underlying ODE system. The equations of the model that were studied in \cite{Fassoni} are
\begin{equation}
\label{início1}
\left\{
\begin{array}{lcl}
\displaystyle	
\frac{d N}{d t} = r_N - \mu_N N - \beta_1 N A - \alpha_N\gamma_N D N,
\vspace{0.2cm}
\\
\displaystyle	
\frac{d A}{d t} = r_A A\left(1-\frac{A}{k_A}\right)-(\mu_A+\epsilon_A)A - \beta_3 NA - \alpha_A\gamma_A D A,
\vspace{0.2cm}
\\
\displaystyle	
\frac{d D}{d t} = \mu - \gamma_A D A - \gamma_N D N - \tau D,
\end{array}
\right.
\end{equation}
where $N$ represents the number of normal cells in a given tissue of the human body, $A$ represents the number of tumor cells in the tissue and $D$ represents the concentration of a chemotherapeutic drug used to treat such a tumor.

Parameter $r_N$ represents a constant influx of new normal cells produced by the tissue stem cells and $\mu_N$ presents the natural mortality of normal cells. A constant influx is considered because the imperative dynamics within a formed tissue is the maintenance of a homeostatic state through the natural replenishment of old and dead cells, see \cite{Simons}.

On the other hand, tumor cells maintain their own growth program \cite{Fedi}. Thus, a density dependent growth is considered for tumor cells. The logistic growth is chosen due to its simplicity. Parameter $\mu_A$ represents the natural mortality of tumor cells, and $\epsilon_A$ represents an extra mortality rate due to apoptosis \cite{Danial}.

Parameters $\beta_1$ and $\beta_3$ encompass the many negative interactions exerted by tumor cells on normal cells and vice-versa, such as competition for nutrients and oxygen. Besides competition, parameter $\beta_3$ encompasses also the effects on normal cells of anti-growth and death signals released by normal cells. In the same way, the parameter $\beta_1$ encompasses also mechanisms developed by tumor cells that damage normal tissue, such as increased local acidity, growth suppression, and release of death signals \cite{Hanahan}.

The third equation of (\ref{início1}) describes the dynamics of chemotherapeutic drug concentration according the following assumptions. The drug has a constant infusion rate $\mu$ and a clearance rate $\tau$. Such constant infusion rate mimics a metronomic dosage, i.e., a near continuous and long-term administration of the drug. The absorption and deactivation of the drug by normal and cancerous cells are described in terms of the law of mass action with rates $\gamma_N$ and $\gamma_A$. Following the log-linear hypothesis \cite{Andre}, it is assumed that the amounts of drug absorbed by normal ($\gamma_N N D$) and cancerous cells ($\gamma_A A D$) kill such cells with rates $\alpha_N$ and $\alpha_A$, respectively. Although many models of cancer treatment do not consider drug absorption explicitly, in \cite{Fassoni}, the authors believe that it is an important fact to be considered, since, this phenomenon contributes to decrease the concentration of drug as time passes.

System (\ref{início1}) is similar to the classical Lotka-Volterra competition model, frequently used in models for tumor growth and population dynamics. The fundamental difference here is the use of a constant flux for normal cells instead of a logistic growth. Such constant flux, also used in other well-known models of cancer \cite{Earn}, removes the symmetry observed in the Lotka-Volterra equations, so that there is no steady state with $N = 0$. Thus, it is impossible to observe the extinction of one of the populations (the normal cells in this case), as opposed to the Lotka-Volterra models. The authors of \cite{Fassoni} claim that this is a realistic result since, roughly speaking, cancer "does not win" by killing all the cells in the tissue, but by reaching a dangerous size that disrupts the proper functioning of the tissue and threatens the health of the individual.

In this work, we are not interested in analyzing the dynamics (stability, asymptotic behavior) of the model, as such study has already been made in \cite{Fassoni}. Our objective is to study the existence and uniqueness of the solution of system (\ref{0riginalEquations}). Such system extends the ODE model \eqref{início1} to a more realistic situation by considering spatial variation of normal and cancer cells and the diffusion of the chemotherapeutic drug through the tissue, with diffusion coefficient $\sigma$ \cite{Anderson}. Further, it is also assumed that the drug influx is restricted to a limited region of the tissue, corresponding to a blood vessel passing transversely in such region. This is mathematically described in the model by the expression $\mu \chi_{\omega} $, where $\chi_{\omega} $ is the characteristic function of the subset $\omega \subset \Omega$. Finally, due to mathematical necessity to simplify the model, we set $\beta_3=0$. This corresponds to a situation where normal cells do not exert negative effects on tumor cells, and is a plausible biological assumption, since there are many tumors that develop resistance to the normal tissue' mechanisms which suppress tumor growth \cite{Hanahan}.

The paper is organized as follows. In Section 2 we present the technical hypothesis and state our main result. In Section 3 we study an auxiliary problem. Using its solution, we prove our main result in Section 4. In Section 5 we present numerical simulations illustrating model behavior.

\section{Technical hypotheses and main result}
Let $\Omega \subset \R^2$ be a domain with boundary $\partial\Omega$, $0\leq T <\infty$, and denote $Q = \Omega\times (0, T)$ and $\Gamma =  \partial \Omega\times (0, T)$.
We will use standard notations for Sobolev spaces, i.e., given $1~\leq~p~\leq~+\infty $ and $k \in\mathbb{N}$, we denote
$$W_{p}^{k}(\Omega)=\left\{ f \\ \in L^{p}(\Omega) : D^{\alpha}f \in L^{p}(\Omega), |\alpha| \leq k \right\}; $$
when $p=2$, as usual we denote $W_{2}^{k}(\Omega) = H^k (\Omega)$;
properties of these spaces can be found for instance in Adams~\cite[Theorem~5.4, p. 97]{Adams}.
Problem~(\ref{0riginalEquations}) will be studied in the standard functional spaces denoted by
\begin{eqnarray*}
  W_{q}^{2,1}(Q) & =& \left\{f\in L^{q}(Q):D^{\alpha}f\in L^{q}(Q), \, \forall
1\leq|\alpha|\leq 2,  f_t \in L^{q}(Q)\right\},
\end{eqnarray*}
\begin{eqnarray*}
  W &=& \left\{f\in L^\infty(Q): f_t \in L^\infty(Q)\right\}
\end{eqnarray*}
and
\begin{eqnarray*}
  L^{p}(0,T;B) &=&  \left\{f:(0,T)\rightarrow B:   \|f(t)\|_{L^{p}(0,T;B)} <+\infty \right\},
\end{eqnarray*}
where $B$ is suitable Banach space, and the norm is given by
$\|f(t)\|_{L^{p}(0,T;B)} = \|\ \|f(t)\|_{B}\ \|_{L^{p}((0,T))}$.
We remark that $L^{p}(Q) = L^{p}((0,T);L^{p}(\Omega))$.
Results concerning these spaces can be found  for instance in Ladyzhenskaya~\cite{Ladyzhenskaya}
and Mikhaylov~\cite{Mikhaylov}.

\vspace{0.1cm}
Next, we state some hypotheses that will be assumed throughout this article.

\subsection{Technical Hypotheses:}
\label{MainHypotheses}

\begin{itemize}

\item[{\bf (i)}] $\Omega\subset\mathbb{R}^2$ is a bounded $C^2$-domain;

\item[{\bf (ii)}] $0< T < \infty$, and $Q=\Omega\times(0,T)$;

\item[{\bf (iii)}] $N_0, A_0 \in L^{\infty}(\Omega)$ and $D_0 \in W^{\frac{3}{2}}_{4}(\Omega)$, satisfying $\frac{\partial D_0}{\partial \eta} (\cdot) =0, \textup{ on } \partial\Omega$;

\item[{\bf (iv)}] $0 \leq D_0 \leq \frac{\mu}{\tau}$ and $N_0, A_0 \ge 0$ a.e. on $\Omega$.

\end{itemize}

\begin{remark}
The constraints imposed in~{\bf (iv)} on the initial conditions  are natural biological requirements.
\end{remark}

\subsection{Main result:}
\begin{theorem}
\label{Teorema1}
Assume that the Technical Hypotheses \ref{MainHypotheses} hold;
then, there exists a unique nonnegative solution $(N,A,D) \in W \times W \times W^{2,1}_4(Q)$ of Problem (\ref{0riginalEquations}). Moreover, $N, A$ and $D$ are functions satisfying
\begin{eqnarray*}
  N \leq ||N_0||_{L^\infty(Q)} + r_N T, \ A \leq C_{\lambda}||A_0||_{L^\infty(\Omega)} \ a.e. \ in \ Q
\end{eqnarray*}
and
\begin{eqnarray*}
  ||N||_W + ||A||_{W} + ||D||_{ W^{2,1}_4(Q)} \leq C,
\end{eqnarray*}
where $C$ is a constant depending on $r_N$, $\mu_N$, $\beta_1$, $\alpha_N$, $\gamma_N$, $C_{\lambda}$, $r_A$, $k_A$, $\mu_A$, $\epsilon_A$, $\alpha_A$, $\gamma_A$, $\mu$, $\tau$, $T$, $\omega$, $||N_0||_{L^\infty (\Omega)}$, $||A_0||_{L^\infty(\Omega)}$ and $||D_0||_{W^{\frac{3}{2}}_4(\Omega)}$.
\end{theorem}

\begin{remark}
The explicit knowledge on how the constant $C$ appearing in the above estimates depends on the given data is important for applications in related control problems.
\end{remark}

\subsection{Known technical results:}

To ease the references, we also state some technical results to be used in this paper.
The first one is sometimes called the Lions-Peetre embedding theorem (see Lions~\cite{Lions}, pp.15);
it is also a particular case of Lemma~3.3, pp.80, in Ladyzhenskaya~\cite{Ladyzhenskaya}:
(obtained by taking $l = 1$ and $r = s = 0$).
\begin{lemma}
\label{icontLp01}
Let $\Omega$ be a domain of $\R^n$ with boundary $\partial \Omega$ satisfying the cone property.
Then, the functional space $W^{2,1}_p(Q)$ is continuously embedded in $u \in L^{q}(Q)$ for $q$ satisfying: {\bf (i)} $1 \leq q \leq \frac{p(n+2)}{n+2-2p}$, if $ p< \frac{n+2}{2}$; {\bf (ii)} $1 \leq q <\infty$, if $p= \frac{n+2}{2}$ and {\bf (iii)} $q=\infty$, if $p>\frac{n+2}{2}$.
\noindent
In particular, for such $q$ and any function $u \in W^{2,1}_p(Q)$ we have that
\begin{eqnarray*}
\label{i.01}
\displaystyle \|u\|_{L^{q}(Q)} \leq C\|u\|_{W^{2,1}_p(Q)},
\end{eqnarray*}
\noindent
with a constant $C$ depending only on $\Omega$, $T$, $p$, $q$, $n$.

In the cases {\bf (ii)}, {\bf (iii)} or in {\bf (i)} when $\displaystyle 1 \leq q < \frac{p(n+2)}{n+2-2p}$, the referred embedding is compact.

\end{lemma}

\vspace{0.1cm}
Next, we consider the following simple parabolic initial-boundary value problem:
\begin{equation}
\label{P_Newmman}
\left\{
\begin{array}{lcl}
\displaystyle	\frac{\partial u}{\partial t} - \sum\limits_{i,j=1}^na_{ij}(x,t)\frac{\partial u^2}{\partial x_ix_j}
+ \sum\limits_{j=1}^na_i(x,t)\frac{\partial u}{\partial x_j} + a(x,t)u=f & \textup{in} & Q,
\\
\displaystyle
\sum\limits_{i=1}^n b_i(x,t)\frac{\partial u}{\partial x_i} + b(x,t)u =0
& \textup{on} & \Gamma,
\\
\displaystyle
u(\cdot,0)= u_0(\cdot) & \textup{in} & \Omega .
 \end{array}
 \right.
\end{equation}

Existence and uniqueness of solutions for this problem is a particular case of Theorem~$9.1$, pp.$341$,
in Ladyzenskaya~\cite{Ladyzhenskaya} for the case of Neumann boundary condition, according to the remarks at the end Chapter IV, section 9, p. 351 in \cite{Ladyzhenskaya}.
In the following, we state this particular result, stressing the dependencies certain norms of the coefficients,  that will be important in our future arguments.
\begin{proposition}
\label{sol. Neumann}
Let $\Omega$ be a bounded domain in $\mathbb{R}^n$, with a $C^{2}$ boundary $\partial \Omega$,
$a_{ij}$ be bounded continuous functions in $Q$,
and $q > 1$. Assume that
\begin{enumerate}

\item $a_{ij} \in C(\bar{Q})$, $i, j=1, \ldots, n$; $[a_{ij}]_{n \times n}$ is a real positive matrix such that for some positive constant $\beta$
we have  $ \sum\limits_{i,j=1}^n a_{ij}(x,t)\xi_i\xi_j\geq \beta|\xi|^2$ for all
$(x,t) \in Q$ and all $\xi \in R^n$,
;

\item $\displaystyle f \in L^p(Q)$;

\item $\displaystyle a_i \in L^r(Q)$ with
either $r =  \max\big(p, n + 2\big)$ if $p \neq n + 2$
or $r =  n + 2 + \varepsilon$, for any  $\varepsilon>0$, if $p = n + 2$;

\item $\displaystyle a \in L^s(Q)$ with
either $s = \max\big(p, (n + 2)/2\big)$ if  $\displaystyle p \neq (n + 2)/2$
or $s = (n + 2)/2 + \varepsilon$, for any $\varepsilon>0$, if $\displaystyle p = (n + 2)/2$.

\item $b_i, b \in C^2 (\bar{\Gamma})$, $i=1, \ldots, n$, and
the coefficients $b_i(x,t)$ satisfy the condition
$\left| \sum\limits_{i=1}^n b_i(x,t)\eta_i(x) \right|\geq \delta >0$
for $a.e.$ in  $\partial\Omega \times (0,T)$,
where $\eta_i(x)$ is the $i^{th}$-component of the unitary outer normal vector to $\partial\Omega$ in $x \in \partial\Omega$;

\item $u_0 \in \ W^{2 - \frac{2}{p}}_p(\Omega)$ with $p\neq 3$ and satisfying the compatibility condition
\\
$\displaystyle \sum\limits_{i=1}^n b_i \frac{\partial u_0}{\partial x_i} + b \ u_0 =0$ on $\partial \Omega$ when $p > 3$.

\end{enumerate}
Then, there exists a unique solution $u \in W^{2,1}_p(Q)$ of Problem~(\ref{P_Newmman});
moreover, there is a positive constant $C_p$  such that
the solution satisfies
\begin{equation}
\label{BasicParabolicEstimate}
\|u\|_{W^{2,1}_p(Q)} \leq C_{p} \left(\|f\|_{L^p(Q)} + \|u_0\|_{W^{2  - \frac{2}{p}}_p(\Omega)}\right).
\end{equation}

Such constant $C_{p}$ depends only on  $\Omega$, $T$, $p$, $r$, $s$, $\beta$, $\delta$
and on the norms $\|b_i\|_{C^2 (\bar{\Gamma})}$, $\|b\|_{C^2 (\bar{\Gamma})}$, $\|a_{ij}\|_{C(\bar{Q})}$,
$\|a_i\|_{L^r(Q)}$ and $\|a\|_{L^s(Q)}$.
Moreover, we may assume that the dependencies of $C_{p}$ on stated the norms are non decreasing.
\end{proposition}

\begin{remark}
The result set out in Proposition \ref{sol. Neumann} can be formulated for the parabolic problem with Dirichlet conditions (see Ladyzenskaya \cite[Theorem 9.1, pp.$341$]{Ladyzhenskaya}). In the problem with Dirichlet condition the compatibility condition in Proposition \ref{sol. Neumann}-($6$) can be replaced by $u_0=0$ on $\partial \Omega$ when $p > 3/2$. This way, all the results in this paper holds if we replaced the Neumann conditions by Dirichlet conditions.
\end{remark}

\section{An auxiliary problem}

In this section we will prove an auxiliary result  to be used in the proof of Theorem~\ref{Teorema1}.
To cope with difficulties with the signs of certain terms during the derivation of the estimates,
we firstly have to consider the following modified problem:
\begin{equation}
\label{P01}
\left\{
\begin{array}{lcl}
\displaystyle	
\frac{\partial \hat{N}}{\partial t} = r_N - \mu_N \hat{N} - \beta_1 \hat{N} \hat{A} - \alpha_N\gamma_N |\hat{D}| \hat{N} , & \textup{in}& Q,
\vspace{0.2cm}
\\
\displaystyle	
\frac{\partial \hat{A}}{\partial t} = r_A\hat{A}\left(1-\frac{\hat{A}}{k_A}\right)-(\mu_A+\epsilon_A)\hat{A} - \alpha_A\gamma_A |\hat{D}| \hat{A}, &\textup{in}& Q,
\vspace{0.2cm}
\\
\displaystyle	
\frac{\partial \hat{D}}{\partial t} = \sigma \Delta\hat{D} + \mu\chi_{\omega} - \gamma\hat{D}\hat{A} - \gamma_N\hat{D}\hat{N} - \tau\hat{D}, &\textup{in}& Q,
\vspace{0.2cm}
\\
\displaystyle \frac{\partial \hat{D}}{\partial \eta} (\cdot) =0, &\textup{on}& \Gamma,
\vspace{0.2cm}
\\
\displaystyle \hat{N}(\cdot,0) = N_{0} (\cdot), \hat{A}(\cdot,0) = A_{0}(\cdot), \hat{D}(\cdot,0) = D_{0} (\cdot), &\textup{in}& \Omega.
\end{array}
\right.
\end{equation}

Now we observe that, since the equation for $\hat{N}$ in this last problem is, for each $x \in \Omega$, an ordinary differential equation
which is linear in $\hat{N}$, we can find an explicit expression for it in terms of $\hat{A}$ and $|\hat{D}|$. However, $\hat{A} $ is, for each $x \in \Omega$, a nonlinear differential equation in $\hat{A}$, and we can determine its explicit expression in terms of $|\hat{D}|$ using Bernoulli's method.
Using these observations and setting $\lambda = r_A - (\mu_A + \epsilon_A)$, we introduce operators $\Lambda: L^{\infty}(Q) \to L^{\infty}(Q)$ and $\Theta: L^{\infty}(Q) \to L^{\infty}(Q)$, defined respectively by
\begin{equation}
\label{P77}
\displaystyle	
\Lambda(\phi)(x,t) = \frac{A_0(x)k_A e^{\lambda t} e^{-\alpha_A\gamma_A \int_{0}^{t}|\phi(\xi,x)| d\xi}}{k_A + A_0(x) r_A \int_{0}^{t} e^{\lambda s} e^{-\alpha_A\gamma_A \int_{0}^{s}|\phi(\xi,x)| d\xi} ds}
\end{equation}
and
\begin{equation}
\label{P7}
\displaystyle	
\Theta(\phi)(x,t) = \frac{N_0(x) + r_N \int_{0}^{t} e^{\mu_N s} e^{\alpha_N\gamma_N\int_{0}^{s} |\phi(x,\xi)| d\xi}e^{\beta_1\int_{0}^{s}\Lambda(\phi)(x, \xi) d\xi} ds}{e^{\mu_N t} e^{\alpha_N\gamma_N\int_{0}^{t}|\phi(x,\xi)| d\xi}e^{\beta_1\int_{0}^{t}\Lambda(\phi)(x, \xi) d\xi}},
\end{equation}

\noindent
where $0 \leq s \leq t \leq T$.

\begin{remark}
\label{obs1}
Thus,  $(\hat{N},\hat{A}, \hat{D})$ is a solution of (\ref{P01}) if, and only if, $\hat{N} = \Theta (\hat{D})$, $\hat{A} = \Lambda(\hat{D})$ and $\hat{D}$ satisfies the following integro-differential system:
\begin{equation}
\label{P3}
\left\{
\begin{array}{lcl}
\displaystyle	
\frac{\partial \hat{D}}{\partial t} = \sigma \Delta\hat{D} + \mu\chi_{\omega} - \gamma\hat{D}\Lambda(\hat{D}) - \gamma_N\hat{D} \Theta(\hat{D}) - \tau\hat{D}, &\textup{in}& Q,
\vspace{0.2cm}
\\
\displaystyle \frac{\partial \hat{D}}{\partial \eta} (\cdot) =0, &\textup{on}& \Gamma,
\vspace{0.2cm}
\\
\displaystyle \hat{D}(\cdot,0) = D_{0} (\cdot), &\textup{in}& \Omega.
\end{array}
\right.
\end{equation}
\end{remark}

\begin{remark}
\label{obs2}
Notice that, to guarantee that $(N, A, D)$, with $D = \hat{D}$, $N = \Theta(\hat{D})$ and $A = \Lambda(\hat{D})$ is also a solution of system~(\ref{0riginalEquations}), it is enough to prove that the solution $\hat{D}$ of Problem~(\ref{P3}) is nonnegative.
\end{remark}

For the Problem \ref{P3}, we have the following existence result:
\begin{proposition}\label{Prop1}
Assuming that the Technical Hypotheses~\ref{MainHypotheses} hold, there exists at least one nonnegative solution
$\hat{D} \in W^{2,1}_4(Q)$ of Problem \eqref{P3}. Moreover, $\hat{D} \leq \frac{\mu}{\tau}$ a.e. in $Q$ and
\begin{eqnarray*}
   ||\hat{D}||_{W^{2,1}_4(Q)} \leq C,
\end{eqnarray*}
where $C$ is a constant depending on $\mu$, $T$, $\omega$ and $||D_0||_{W^{\frac{3}{2}}_4(\Omega)}$.
\end{proposition}

\begin{lemma}
\label{base}
Let $f:(0,T) \to \mathbb{R}$ differentiable such that $f(t) > 0$ and $f'(t) \ge 0$. If $g(t) = \frac{\int_{0}^{t} f(x) dx}{f(t)}$, then $g(t) \leq T$, for all $t \in (0, T)$.
\end{lemma}

\noindent
{\bf Proof:}
Since $f$ is continuous in $(0, T)$, it follows that
\begin{eqnarray*}
  g'(t) &=& \frac{f(t)^2 - f'(t) \int_{0}^{t}f(x)dx}{f(t)^2} \\
        &=& 1 - \frac{f'(t)}{f(t)} g(t).
\end{eqnarray*}

As $f(t) > 0$ we have $g(t) \ge 0$ and using the fact that $f'(t) \ge 0$ we obtain $\frac{f'(t)}{f(t)} g(t) \ge 0$. Therefore, $g'(t) \leq 1$, which suggests $g(t) \leq t$, for all $t \in (0, T)$. Thus, $g(t) \leq T$, as intended.
\hfill$\Box$

Since in the proof of existence of solutions of (\ref{P3})
the expression of $\Lambda$ and $\Theta$ will play important roles,
we state some of their properties in the following:
\begin{lemma}
\label{PropertiesEtcFirst}
If $N_0, A_0 \in L^\infty(\Omega)$ and $C_{\lambda} = \max\{1, e^{\lambda T}\}$, then for any $\phi, \phi_1, \phi_2 \in L^\infty(Q)$ and for almost every $(x,t) \in Q$, there holds
\[
\begin{array}{ll}
\mbox{\bf (i)} & 0 \leq \Theta(\phi)(x,t) \leq ||N_0||_{L^\infty(\Omega)} + r_N T;
\vspace{0.2cm}
\\
\mbox{\bf (ii)} & 0 \leq \Lambda(\phi)(x,t) \leq C_{\lambda}||A_0||_{L^\infty(\Omega)};
\vspace{0.2cm}
\\
\mbox{\bf (iii)} & \|\Lambda(\phi_1) - \Lambda(\phi_2) \|_{L^\infty {(Q)}}  \leq C_1 \|\phi_1 - \phi_2\|_{L^\infty(Q)}, \\
& where \ C_1 \ is \ a \ constant \ depending \ on \ r_A, k_A, \alpha_A, \gamma_A, C_{\lambda}, T \ and \ ||A_0||_{L^\infty(\Omega)};
\vspace{0.2cm}
\\
\mbox{\bf (iv)} & \|\Theta(\phi_1) - \Theta(\phi_2) \|_{L^\infty {(Q)}}  \leq C_2 \|\phi_1 - \phi_2\|_{L^\infty(Q)}, \\
& where \ C_2 \ is \ a \ constant \ depending \ on \ r_N, \mu_N, \beta_1,  \alpha_N, \gamma_N, C_{\lambda}, C_1, T , \\
&  ||\phi_1||_{L^\infty(Q)}, ||\phi_2||_{L^\infty(Q)}, ||N_0||_{L^\infty(\Omega)} \ and \ ||A_0||_{L^\infty(\Omega)}.
\end{array}
\]
\end{lemma}

\noindent {\bf Proof (i) and (ii):}
By the expressions (\ref{P77}) and (\ref{P7}) it is immediate that $\Lambda(\phi)(x,t), \Theta(\phi)(x,t) \ge 0$. To prove that $\Theta(\phi)(x,t) \leq ||N_0||_{L^\infty(\Omega)} + r_N T$, we observe that
$$
\begin{array}{rcl}
\displaystyle
\Theta(\phi)(x,t) =  \frac{N_0(x) + r_N \int_{0}^{t} e^{\mu_N s} e^{\alpha_N\gamma_N\int_{0}^{s} |\phi(x,\xi)| d\xi}e^{\beta_1\int_{0}^{s}\Lambda(\phi)(x, \xi) d\xi} ds}{e^{\mu_N t} e^{\alpha_N\gamma_N\int_{0}^{t}|\phi(x,\xi)| d\xi}e^{\beta_1\int_{0}^{t}\Lambda(\phi)(x, \xi) d\xi}}  \\
\displaystyle	
\\
\leq N_0(x) + r_N \frac{\int_{0}^{t} e^{\mu_N s} e^{\alpha_N\gamma_N\int_{0}^{s} |\phi(x,\xi)| d\xi}e^{\beta_1\int_{0}^{s}\Lambda(\phi)(x, \xi) d\xi} ds}{e^{\mu_N t} e^{\alpha_N\gamma_N\int_{0}^{t}|\phi(x,\xi)| d\xi}e^{\beta_1\int_{0}^{t}\Lambda(\phi)(x, \xi) d\xi}}.
\end{array}
$$

Fixed $x \in \Omega$, we define
\begin{equation*}
\displaystyle
g(x,t) = \frac{\int_{0}^{t} e^{\mu_N s} e^{\alpha_N\gamma_N\int_{0}^{s} |\phi(x,\xi)| d\xi}e^{\beta_1\int_{0}^{s}\Lambda(\phi)(x, \xi) d\xi} ds}{e^{\mu_N t} e^{\alpha_N\gamma_N\int_{0}^{t}|\phi(x,\xi)| d\xi}e^{\beta_1\int_{0}^{t}\Lambda(\phi)(x, \xi) d\xi}},
\end{equation*}
and using the Lemma \ref{base} with $f(x,t) = e^{\mu_N t} e^{\alpha_N\gamma_N\int_{0}^{t}|\phi(x,\xi)| d\xi}e^{\beta_1\int_{0}^{t}\Lambda(\phi)(x, \xi) d\xi}$, it follows that
\begin{eqnarray*}
\Theta(\phi, \varphi)(x,t) \leq N_0(x) + r_N T \\
\leq ||N_0||_{L^\infty(\Omega)} + r_N T.
\end{eqnarray*}

To prove that $\Lambda(\phi)(x,t) \leq C_{\lambda}||A_0||_{L^\infty(\Omega)}$, note that
$$
\begin{array}{rcl}
\displaystyle	
\Lambda(\phi)(x,t) = \frac{A_0(x)k_A e^{\lambda t} e^{-\alpha_A \gamma_A \int_{0}^{t}|\phi(\xi,x)| d\xi}}{k_A + A_0(x) r_A \int_{0}^{t} e^{\lambda s} e^{-\alpha_A \gamma_A \int_{0}^{s}|\phi(\xi,x)| d\xi} ds} \\
\displaystyle	
\\
\leq \frac{1}{k_A} A_0(x)k_A e^{\lambda t}e^{-\alpha_A \gamma_A \int^{t}_{0}|\phi(x,\xi)|d\xi} \\
\\
\displaystyle	
\leq A_0(x) e^{\lambda t} \leq C_{\lambda} A_0(x) \leq C_{\lambda} ||A_0||_{L^{\infty}(\Omega)}.
\end{array}
$$

{\bf Proof (iii):} We firstly need to observe that, due to the mean value inequality, given any $z_1, z_2 \in \R$,
there is $\theta = \theta(z_1, z_2)$ such that $e^{z_2} - e^{z_1} = e^{(1-\theta) z_1 + \theta z_2 } (z_2 - z_1)$;
in particular, for any $z_1, z_2 \leq 0$ we also have $(1-\theta) z_1 + \theta z_2  \leq 0$ and thus
\begin{equation}
\label{AlgebraicExponentialInequality}
|e^{z_2} - e^{z_1} | \leq |z_2 - z_1|, \quad \forall z_1, z_2 \leq 0 .
\end{equation}

Secondly, we note that by the inequality (\ref{AlgebraicExponentialInequality}) and by $\phi_i \in L^\infty(Q)$, $i = 1,2$, we obtain
\begin{equation}\label{i1}
\begin{array}{rcl}
\big|e^{-\alpha_A\gamma_A\int_{0}^{t}|\phi_1(x,\xi)| d\xi} - e^{-\alpha_A\gamma_A\int_{0}^{t}|\phi_2(x,\xi)| d\xi}\big| & \leq & \big|-\alpha_A\gamma_A\int_{0}^{t}(|\phi_1(x,\xi)| - |\phi_2(x,\xi)|) d\xi\big| \\
\\
&\leq & \alpha_A\gamma_A T ||\phi_1 - \phi_2||_{L^{\infty}(Q)}.
\end{array}
\end{equation}

Thirdly, we observe that
\begin{eqnarray*}
\big|e^{-\alpha_A\gamma_A\int_{0}^{t}|\phi_1(x,\xi)| d\xi}\int_{0}^{t}e^{\lambda s} e^{-\alpha_A\gamma_A \int_{0}^{s}|\phi_2(x,\xi)| d\xi}ds &-& \\
\displaystyle
e^{-\alpha_A\gamma_A\int_{0}^{t}|\phi_2(x,\xi)| d\xi}\int_{0}^{t}e^{\lambda s} e^{-\alpha_A\gamma_A \int_{0}^{s}|\phi_2(x,\xi)| d\xi}ds\big| &\leq& \\
\displaystyle
\big|e^{-\alpha_A\gamma_A\int_{0}^{t}|\phi_1(x,\xi)| d\xi} - e^{-\alpha_A\gamma_A\int_{0}^{t}|\phi_2(x,\xi)| d\xi}\big| \int_{0}^{t}e^{\lambda s} e^{-\alpha_A\gamma_A \int_{0}^{s}|\phi_2(x,\xi)| d\xi} ds &+& \\
\displaystyle
e^{-\alpha_A\gamma_A\int_{0}^{t}|\phi_2(x,\xi)| d\xi} \int_{0}^{t} e^{\lambda s} \big|e^{-\alpha_A\gamma_A \int_{0}^{s}|\phi_1(x,\xi)| d\xi} - e^{-\alpha_A\gamma_A \int_{0}^{s}|\phi_2(x,\xi)| d\xi}\big|ds.
\end{eqnarray*}

How $e^{\lambda T} \leq C_\lambda$ and $e^{-\alpha_A\gamma_A||\phi_2||_{L^\infty(Q)}} \leq 1$, and using study analogous to that done in (\ref{i1}), we obtain that
\begin{equation}
\label{ineq1x}
\begin{array}{rcl}
\displaystyle
\big|e^{-\alpha_A\gamma_A\int_{0}^{t}|\phi_1(x,\xi)| d\xi}\int_{0}^{t}e^{\lambda s} e^{-\alpha_A\gamma_A \int_{0}^{s}|\phi_2(x,\xi)| d\xi}ds &-& \\
\\
\displaystyle
e^{-\alpha_A\gamma_A\int_{0}^{t}|\phi_2(x,\xi)| d\xi}\int_{0}^{t}e^{\lambda s} e^{-\alpha_A\gamma_A \int_{0}^{s}|\phi_1(x,\xi)| d\xi}ds\big| &\leq& \\
\\
\displaystyle
2 \alpha_A\gamma_A C_{\lambda} T^2 ||\phi_1 - \phi_2||_{L^{\infty}(Q)}.
\end{array}
\end{equation}

Finally, the expression in (\ref{P7}) suggests
\begin{eqnarray*}
|\Lambda(\phi_1)(x,t) - \Lambda(\phi_2)(x,t)| & \leq & \\
A_0(x) e^{\lambda t}\bigg|e^{-\alpha_A\gamma_A\int_{0}^{t}|\phi_1(x,\xi)| d\xi} - e^{-\alpha_A\gamma_A\int_{0}^{t}|\phi_2(x,\xi)| d\xi}\bigg| &+& \\
\frac{1}{k_A}{A_0(x)}^2 r_A e^{\lambda t}\bigg|e^{-\alpha_A\gamma_A\int_{0}^{t}|\phi_1(x,\xi)| d\xi}\int_{0}^{t}e^{\lambda s} e^{-\alpha_A\gamma_A \int_{0}^{s}|\phi_2(x,\xi)| d\xi}ds &-& \\
e^{-\alpha_A\gamma_A\int_{0}^{t}|\phi_2(x,\xi)| d\xi}\int_{0}^{t}e^{\lambda s} e^{-\alpha_A\gamma_A \int_{0}^{s}|\phi_1(x,\xi)| d\xi}ds\bigg|,
\end{eqnarray*}
and using the estimates obtained in (\ref{i1}) and (\ref{ineq1x}) and making the possible simplifications, we obtain
$$
\begin{array}{rcl}
|\Lambda(\phi_1)(x,t) - \Lambda(\phi_2)(x,t)| & \leq & \\
\\
||A_0||_{L^{\infty}(\Omega)} C_{\lambda} \alpha_A\gamma_A T ||\phi_1 - \phi_2||_{L^{\infty}(Q)} &+& \\
\\
\frac{2}{k_A} ||A_0||_{L^{\infty}(\Omega)}^2 r_A {C_{\lambda}}^2 \alpha_A\gamma_A T^2 ||\phi_1 - \phi_2||_{L^{\infty}(Q)},
\end{array}
$$
for almost everything $(x,t) \in Q$, i.e.,
\begin{equation}
\label{eqc1}
 ||\Lambda(\phi_1) - \Lambda(\phi_2)||_{L^{\infty}(Q)} \leq  C_1 \ ||\phi_1 - \phi_2||_{L^{\infty}(Q)}.
\end{equation}

{\bf Proof (iv):} First, note that
\begin{eqnarray*}
\big| e^{\alpha_N\gamma_N\int_{0}^{t}|\phi_2(x,\xi)| d\xi}e^{\beta_1\int_{0}^{t}\Lambda(\phi_2)(x, \xi) d\xi} - e^{\alpha_N\gamma_N\int_{0}^{t}|\phi_1(x,\xi)| d\xi}e^{\beta_1\int_{0}^{t}\Lambda(\phi_1)(x, \xi) d\xi}\big| &\leq& \\
e^{\alpha_N\gamma_N\int_{0}^{t}|\phi_2(x,\xi)| d\xi}\big|e^{\beta_1\int_{0}^{t}\Lambda(\phi_2)(x, \xi) d\xi} - e^{\beta_1\int_{0}^{t}\Lambda(\phi_1)(x, \xi) d\xi} \big| &+& \\
e^{\beta_1\int_{0}^{t}\Lambda(\phi_1)(x, \xi) d\xi}\big|e^{\alpha_N\gamma_N\int_{0}^{t}|\phi_2(x,\xi)| d\xi} - e^{\alpha_N\gamma_N\int_{0}^{t}|\phi_1(x,\xi)| d\xi} \big|,
\end{eqnarray*}
and by the inequality (\ref{AlgebraicExponentialInequality}) and by $\Lambda(\phi_i), \phi_i \in L^\infty(Q)$, $i = 1,2$, we obtain
\begin{equation}
\label{i12}
\begin{array}{rcl}
\displaystyle
\big| e^{\alpha_N\gamma_N\int_{0}^{t}|\phi_2(x,\xi)| d\xi}e^{\beta_1\int_{0}^{t}\Lambda(\phi_2)(x, \xi) d\xi} - e^{\alpha_N\gamma_N\int_{0}^{t}|\phi_1(x,\xi)| d\xi}e^{\beta_1\int_{0}^{t}\Lambda(\phi_1)(x, \xi) d\xi}\big| &\leq& \\
\\
\displaystyle
e^{\alpha_N\gamma_N\int_{0}^{t}|\phi_2(x,\xi)| d\xi} \beta_1 T ||\Lambda(\phi_1) - \Lambda(\phi_2)||_{L^\infty(Q)} &+& \\
\\
\displaystyle
e^{\beta_1\int_{0}^{t}\Lambda(\phi_1)(x, \xi) d\xi} \alpha_N \gamma_N T ||\phi_1 - \phi_2||_{L^\infty(Q)}.
\end{array}
\end{equation}

Since
\begin{eqnarray*}
\bigg|e^{\alpha_N\gamma_N\int_{0}^{t}|\phi_2(x,\xi)| d\xi}e^{\beta_1\int_{0}^{t}\Lambda(\phi_2)(x, \xi) d\xi}\int_{0}^{t}e^{\mu_Ns}e^{\alpha_N\gamma_N\int_{0}^{s}|\phi_1(x,\xi)| d\xi}e^{\beta_1\int_{0}^{s}\Lambda(\phi_1)(x, \xi) d\xi} ds &-& \\
\\
\displaystyle
e^{\alpha_N\gamma_N\int_{0}^{t}|\phi_1(x,\xi)| d\xi}e^{\beta_1\int_{0}^{t}\Lambda(\phi_1)(x, \xi) d\xi}\int_{0}^{t}e^{\mu_Ns}e^{\alpha_N\gamma_N\int_{0}^{s}|\phi_2(x,\xi)| d\xi}e^{\beta_1\int_{0}^{s}\Lambda(\phi_2)(x, \xi) d\xi} ds\bigg| &\leq& \\
\\
\displaystyle
e^{\alpha_N\gamma_N\int_{0}^{t}|\phi_2(x,\xi)| d\xi}e^{\beta_1\int_{0}^{t}\Lambda(\phi_2)(x, \xi) d\xi} &\times& \\
\\
\displaystyle
\int_{0}^{t}e^{\mu_Ns}\bigg|e^{\alpha_N\gamma_N\int_{0}^{s}|\phi_1(x,\xi)| d\xi}e^{\beta_1\int_{0}^{s}\Lambda(\phi_1)(x, \xi) d\xi} - e^{\alpha_N\gamma_N\int_{0}^{s}|\phi_2(x,\xi)| d\xi}e^{\beta_1\int_{0}^{s}\Lambda(\phi_2)(x, \xi) d\xi}\bigg| ds &+& \\
\\
\displaystyle
\bigg|e^{\alpha_N\gamma_N\int_{0}^{s}|\phi_2(x,\xi)| d\xi}e^{\beta_1\int_{0}^{s}\Lambda(\phi_2)(x, \xi) d\xi} - e^{\alpha_N\gamma_N\int_{0}^{s}|\phi_1(x,\xi)| d\xi}e^{\beta_1\int_{0}^{s}\Lambda(\phi_1)(x, \xi) d\xi}\bigg| &\times& \\
\\
\displaystyle
\int_{0}^{t}e^{\mu_Ns}e^{\alpha_N\gamma_N\int_{0}^{s}|\phi_2(x,\xi)| d\xi}e^{\beta_1\int_{0}^{s}\Lambda(\phi_2(x, \xi) d\xi},
\end{eqnarray*}
doing $||\phi||_{L^\infty(Q)} = \max\{||\phi_1||_{L^\infty(Q)},  ||\phi_2||_{L^\infty(Q)}\}$ and study analogous to that done in (\ref{i12}), guarantees us
\begin{equation}
\label{ineq21}
\begin{array}{rcl}
\displaystyle
\bigg|e^{\alpha_N\gamma_N\int_{0}^{t}|\phi_2(x,\xi)| d\xi}e^{\beta_1\int_{0}^{t}\Lambda(\phi_2)(x, \xi) d\xi}\int_{0}^{t}e^{\mu_Ns}e^{\alpha_N\gamma_N\int_{0}^{s}|\phi_1(x,\xi)| d\xi}e^{\beta_1\int_{0}^{s}\Lambda(\phi_1)(x, \xi) d\xi} ds &-& \\
\\
\displaystyle
e^{\alpha_N\gamma_N\int_{0}^{t}|\phi_1(x,\xi)| d\xi}e^{\beta_1\int_{0}^{t}\Lambda(\phi_1)(x, \xi) d\xi}\int_{0}^{t}e^{\mu_Ns}e^{\alpha_N\gamma_N\int_{0}^{s}|\phi_2(x,\xi)| d\xi}e^{\beta_1\int_{0}^{s}\Lambda(\phi_2)(x, \xi) d\xi} ds\bigg| &\leq& \\
\displaystyle
\\
e^{\alpha_N\gamma_N\int_{0}^{t}|\phi_2(x,\xi)| d\xi}e^{\beta_1 \int_{0}^{t}\Lambda(\phi_2)(x, \xi) d\xi}e^{\mu T}  e^{\alpha_N\gamma_N T ||\phi||_{L^{\infty}(Q)}} \beta_1 T^2 &\times& \\
\displaystyle
\\
||\Lambda(\phi_1) - \Lambda(\phi_2)||_{L^\infty(Q)} &+& \\
\\
\displaystyle
e^{\alpha_N\gamma_N\int_{0}^{t}|\phi_2(x,\xi)| d\xi}e^{\beta_1\int_{0}^{t}\Lambda(\phi_2)(x, t) d\xi}e^{\mu T}  e^{\beta_1 T C_{\lambda} ||A_0||_{L^\infty(\Omega)}}\alpha_N\gamma_N T^2 &\times& \\
\\
\displaystyle
||\phi_1 - \phi_2||_{L^\infty(Q)} &+& \\
\\
\displaystyle
e^{\alpha_N\gamma_N\int_{0}^{t}|\phi_2(x,\xi)| d\xi} \beta_1 T^2 ||\Lambda(\phi_1) - \Lambda(\phi_2)||_{L^\infty(Q)} &\times& \\
\\
\displaystyle
e^{\mu_N T}e^{\alpha_N\gamma_N T ||\phi||_{L^{\infty}(Q)}}e^{\beta_1 T C_{\lambda} ||A_0||_{L^\infty(\Omega)}} &+& \\
\\
\displaystyle
e^{\beta_1\int_{0}^{t}\Lambda(\phi_1)(x, \xi) d\xi} \alpha_N \gamma_N T^2 ||\phi_1 - \phi_2||_{L^\infty(Q)} &\times& \\
\\
\displaystyle
e^{\mu_N T}e^{\alpha_N\gamma_N T ||\phi||_{L^{\infty}(Q)}}e^{\beta_1 T C_{\lambda} ||A_0||_{L^\infty(\Omega)}}.
\end{array}
\end{equation}

Finally, the expression in (\ref{P7}) suggests
\begin{eqnarray*}
|\Theta(\phi_1)(x,t) - \Theta(\phi_2)(x,t)| & \leq & \\
\frac{1}{e^{\alpha_N\gamma_N\int_{0}^{t}|\phi_1(x,\xi)| d\xi}e^{\beta_1\int_{0}^{t}\Lambda(\phi_1)(x, \xi) d\xi} e^{\alpha_N\gamma_N\int_{0}^{t}|\phi_2(x,\xi)| d\xi}e^{\beta_1\int_{0}^{t}\Lambda(\phi_2)(x, \xi) d\xi}}
&\times& \\
\bigg(N_0(x)\big|e^{\alpha_N\gamma_N\int_{0}^{t}|\phi_2(x,\xi)| d\xi}e^{\beta_1\int_{0}^{t}\Lambda(\phi_2)(x, \xi) d\xi} - e^{\alpha_N\gamma_N\int_{0}^{t}|\phi_1(x,\xi)| d\xi}e^{\beta_1\int_{0}^{t}\Lambda(\phi_1)(x, \xi) d\xi} \big| &+& \\
r_N\bigg|e^{\alpha_N\gamma_N\int_{0}^{t}|\phi_2(x,\xi)| d\xi}e^{\beta_1\int_{0}^{t}\Lambda(\phi_2)(x, \xi) d\xi}\int_{0}^{t}e^{\mu_Ns}e^{\alpha_N\gamma_N\int_{0}^{s}|\phi_1(x,\xi)| d\xi}e^{\beta_1\int_{0}^{s}\Lambda(\phi_1)(x, \xi) d\xi} ds &-& \\
e^{\alpha_N\gamma_N\int_{0}^{t}|\phi_1(x,\xi)| d\xi}e^{\beta_1\int_{0}^{t}\Lambda(\phi_1)(x, \xi) d\xi}\int_{0}^{t}e^{\mu_Ns}e^{\alpha_N\gamma_N\int_{0}^{s}|\phi_2(x,\xi)| d\xi}e^{\beta_1\int_{0}^{s}\Lambda(\phi_2)(x, \xi) d\xi}\bigg|\bigg).
\end{eqnarray*}
and using the estimates obtained in (\ref{eqc1}), (\ref{i12}) and (\ref{ineq21}) and making the possible simplifications, we obtain
$$
\begin{array}{rcl}
\displaystyle
\displaystyle
|\Theta(\phi_1)(x,t) - \Theta(\phi_2)(x,t)| & \leq & \\
\\
\displaystyle
||N_0||_{L^{\infty}(\Omega)} e^{\alpha_N\gamma_N T ||\phi||_{L^{\infty}(Q)}} \beta_1 T^2 C_1 ||\phi_1 - \phi_2||_{L^{\infty}(Q)} &+& \\
\\
\displaystyle
||N_0||_{L^{\infty}(\Omega)} e^{\beta_1 T C_{\lambda} ||A_0||_{L^\infty(\Omega)}}\alpha_N\gamma_N T^2 ||\phi_1 - \phi_2||_{L^{\infty}(Q)} &+& \\
\\
\displaystyle
r_N e^{\mu_N T} e^{\alpha_N\gamma_N T ||\phi||_{L^{\infty}(Q)}} e^{\beta_1 T C_{\lambda} ||A_0||_{L^\infty(\Omega)}} \beta_1 T^2 C_1 ||\phi_1 - \phi_2||_{L^{\infty}(Q)}  &+& \\
\\
\displaystyle
r_N e^{\mu_N T} e^{\alpha_N\gamma_N T ||\phi||_{L^{\infty}(Q)}} e^{\beta_1 T C_{\lambda} ||A_0||_{L^\infty(\Omega)}}  \alpha_N\gamma_N T^2 ||\phi_1 - \phi_2||_{L^{\infty}(Q)}
\end{array}
$$
for almost everything $(x,t) \in Q$, i.e.,
\begin{eqnarray*}
\displaystyle
||\Theta(\phi_1) - \Theta(\phi_2)||_{L^\infty(Q)} \leq C_2 ||\phi_1 -  \phi_2||_{L^\infty(Q)}.
\end{eqnarray*}
\hfill$\Box$

\subsection{Proof of Proposition \ref{Prop1}}

To not overburden the notation, in this subsection we denote $D$ as a generic solution of the equations that follows.

To get a solution of problem \eqref{P3}, we will apply the Leray-Schauder fixed point theorem to the mapping
$\Psi$ defined as follows:
\begin{equation}
\label{oper}
\begin{array}{rccl}
\Psi: & [0,1]\times L^\infty(Q)  & \rightarrow & L^\infty(Q)
\\
&(l, \phi) & \mapsto &	D,
\end{array}
\end{equation}

\noindent
where $D$ is the unique solution of
\begin{equation}
\label{P6}
\left\{
\begin{array}{lcl}
\displaystyle	
\frac{\partial D}{\partial t} = \sigma \Delta D + \mu\chi_{\omega} - l\gamma D\Lambda(\phi) - l\gamma_N D \Theta(\phi) - \tau D, &\textup{in}& Q,
\vspace{0.2cm}
\\
\displaystyle \frac{\partial D}{\partial \eta} (\cdot) =0, &\textup{on}& \Gamma,
\vspace{0.2cm}
\\
\displaystyle D(\cdot,0) = D_{0} (\cdot), &\textup{in}& \Omega,
\end{array}
\right.
\end{equation}
with $\Lambda(\phi)$ and $\Theta(\phi)$ given by (\ref{P77}) and (\ref{P7}), respectively.

To apply such theorem we present next a sequence of lemmas:

\begin{lemma}
\label{lemaAi}
Suppose $N_0, A_0 \in L^\infty(\Omega)$ and $D_0 \in W^{\frac{3}{2}}_4(\Omega)$. Then the mapping $\Psi :[0,1] \times L^\infty(Q) \rightarrow L^\infty(Q)$ is well defined.
\end{lemma}

\noindent{\bf Proof:}
We affirm that the coefficients of the Problem \ref{P6} satisfy the hypotheses of the Proposition \ref{sol. Neumann}. For example, it is immediate that $- l \gamma \Lambda(\phi) - l \gamma_N \Theta(\phi) - \tau \in L^4(Q)$, because by Lemma \ref{PropertiesEtcFirst}, $\Lambda(\phi), \Theta(\phi) \in L^\infty(Q)$. Thus, we conclude that there is a unique solution $D \in W^{2,1}_4(Q)$ of problem $\ref{P6}$. Moreover, $D$ satisfies the following estimate:
\begin{equation}\label{AA}
\begin{array}{cccc}
  ||D||_{W^{2,1}_4(Q)} \leq C_p \bigg( ||\mu \chi_{\omega}||_{L^4(Q)} + ||D_0||_{W^{\frac{3}{2}}_4(\Omega)}\bigg) \\
  \leq C_p \bigg( \mu |\omega|^{\frac{1}{4}} T^{\frac{1}{4}} + ||D_0||_{W^{\frac{3}{2}}_4(\Omega)} \bigg).
\end{array}
\end{equation}

Finally, from Lemma \ref{icontLp01}, we have $W^{2,1}_4(Q) \hookrightarrow L^\infty(Q)$, and we conclude that the operator $\Psi$ in well defined.
\hfill$\Box$
\\

\begin{lemma}
\label{nãonegativa1}
Suppose $D$ is a solution of (\ref{P6}) and $0 \leq D_0 \leq \frac{\mu}{\tau}$ a.e. in $\Omega$, then $0 \leq D \leq  \frac{\mu}{\tau}$ a.e. in $Q$.
\end{lemma}

\noindent {\bf Proof:}
Multiplying the first equation in $(\ref{P6})$ by $D^-$ and integrating into $\Omega$, we get
\begin{eqnarray*}
   \frac{1}{2} \frac{d}{dt} \int_{\Omega} (D^{-})^2 \ dx  = -\sigma  \int_{\Omega} |\nabla D^{-}|^2 \ dx - \mu \int_{\omega} D^{-} \ dx \\
   - l\gamma \int_{\Omega} \Lambda(\phi) (D^{-})^2 \ dx - l\gamma_N \int_{\Omega} \Theta(\phi) (D^{-})^2 \ dx - \tau \int_{\Omega} (D^{-})^2 \ dx.
\end{eqnarray*}

Thus,
\begin{eqnarray*}
\frac{d}{dt} \int_{\Omega} (D^{-})^2 dx \leq 0,
\end{eqnarray*}
and using Gronwall's inequality and the fact that $D_0 \ge 0$ a.e. in $\Omega$, we obtain
\begin{eqnarray*}
\int_{\Omega} (D^{-})^2 dx \leq \int_{\Omega}
({D_0}^-)^2 dx = 0,
\end{eqnarray*}
that is, $||D^-(\cdot, t)||_{L^2(\Omega)} = 0$ for all $t \in (0,T)$, where we conclude that $D^{-} = 0$ a.e. in $Q$ and therefore $D
\ge 0$ a.e. in $Q$.

Now, we observe that the first equation in (\ref{P6}) can be rewritten as
\begin{eqnarray*}
\frac{\partial}{\partial t} \big(D - \frac{\mu}{\tau}\big) = \sigma \Delta \big(D - \frac{\mu}{\tau}\big)- l\gamma \Lambda(\phi) D
- l\gamma_N \Theta(\phi) D - \tau\big(D - \frac{\mu\chi_{\omega}}{\tau} \big).
\end{eqnarray*}

Multiplying by $(D - \frac{\mu}{\tau})^+$ and integrating in $\Omega$, we obtain
\begin{eqnarray*}
\frac{1}{2} \frac{d}{dt} \int_{\Omega} \big(\big(D - \frac{\mu}{\tau}\big)^{+}\big)^2 \ dx = - \sigma \int_{\Omega} \big|\nabla \big(D - \frac{\mu}{\tau}\big)^{+} \big|^2 \ dx \\
- l\gamma \int_{\Omega} \Lambda(\phi) D \big(D - \frac{\mu}{\tau}\big)^{+} \ dx - l\gamma_N \int_{\Omega} \Theta(\phi) D \big(D - \frac{\mu}{\tau}\big)^{+} \ dx \\
- \tau \int_{\omega} \big(\big(D - \frac{\mu}{\tau} \big)^{+}\big)^2 \ dx - \tau \int_{\Omega \backslash \omega} D \big(D - \frac{\mu}{\tau} \big)^{+} \ dx,
\end{eqnarray*}
that is,
\begin{eqnarray*}
  \frac{d}{dt} \int_{\Omega} \big(\big(D - \frac{\mu}{\tau}\big)^{+}\big)^2 \ dx \leq 0.
\end{eqnarray*}

Thus, using Gronwall's inequality and the fact that $D_0 \leq \frac{\mu}{\tau}$ a.e. in $\Omega$, it follows that
\begin{eqnarray*}
\int_{\Omega} \big(\big({D} -\frac{\mu}{\tau}\big)^{+}\big)^2 \ dx &\leq& \int_{\Omega} \big(\big(D_0 - \frac{\mu}{\tau} \big)^{+}\big)^2 \ dx = 0,
\end{eqnarray*}
that is, $||\big(D(\cdot, t) -\frac{\mu}{\tau}\big)^{+}||_{L^2(\Omega)} = 0$ for all $t \in (0,T)$, and therefore $\big(D -\frac{\mu}{\tau}\big)^{+} = 0$ a.e. in $Q$, and we conclude that $D \leq \frac{\mu}{\tau}$ a.e. in $Q$.

\hfill$\Box$
\\

\begin{lemma}
\label{con1}
For each fixed $l \in [0,1]$, the mapping $\Psi(l, \cdot): L^\infty(Q) \rightarrow L^\infty(Q)$ is compact, i.e.,
it is continuous and maps bounded sets into relatively compacts sets.
\end{lemma}

\noindent{\bf Proof:}
The functions $\Psi(l, \phi_1)= D_1$ and $\Psi(l, \phi_2) = D_2$ satisfy the system
\begin{equation*}
\label{Aij}
\left\{
\begin{array}{lcl}
\displaystyle	
\frac{\partial D_i}{\partial t} = \sigma \Delta D_i + \mu\chi_{\omega} - l\gamma D_i\Lambda(\phi_1) - l\gamma_N D_i \Theta(\phi_i) - \tau D_i, &\textup{in}& Q,
\vspace{0.2cm}
\\
\displaystyle \frac{\partial D_i}{\partial \eta} (\cdot) =0, &\textup{on}& \Gamma,
\vspace{0.2cm}
\\
\displaystyle D_i(\cdot,0) = D_{0} (\cdot), &\textup{in}& \Omega,
\end{array}
\right.
\end{equation*}
with $i=1,2$; letting $\tilde{D} = D_1 - D_2$, we have
\begin{equation}
\label{Ai1}
\left\{
\begin{array}{lcl}
\displaystyle	
\frac{\partial \tilde{D}}{\partial t} - \sigma \Delta \tilde{D} + l\gamma \tilde{D} \Lambda(\phi_2) + l\gamma_N \tilde{D} \Theta(\phi_2) + \tau \tilde{D} =
\\
\displaystyle  - l\gamma D_1 (\Lambda(\phi_1) - \Lambda(\phi_2)) - l\gamma_N D_1(\Theta(\phi_1) - \Theta(\phi_2)), &\textup{in}& Q,
\vspace{0.2cm}
\\
\displaystyle \frac{\partial \tilde{D}}{\partial \eta} (\cdot) =0, &\textup{on}& \Gamma,
\vspace{0.2cm}
\\
\displaystyle \tilde{D}(\cdot,0) = \tilde{D}_0 (\cdot) = 0, &\textup{in}& \Omega.
\end{array}
\right.
\end{equation}

Using the Proposition \ref{sol. Neumann} and the fact that $L^{\infty}(Q) \hookrightarrow L^4(Q)$ and $D_1 \leq \frac{\mu}{\tau}$, we get
\begin{eqnarray*}
 ||\tilde{D}||_{W_{4}^{2, 1}(Q)} \leq C_p  ||- l\gamma D_1 (\Lambda(\phi_1) - \Lambda(\phi_2))
- l\gamma_N D_1(\Theta(D_1) - \Theta(D_2)) ||_{L^4(Q)} \\
\leq \bar{C}_p ||- l\gamma D_1 (\Lambda(\phi_1) - \Lambda(\phi_2)) - l\gamma_N D_1(\Theta(\phi_1) - \Theta(\phi_2)) ||_{L^\infty(Q)} \\
\leq \bar{C}_p \gamma \frac{\mu}{\tau} ||\Lambda(\phi_1) - \Lambda(\phi_2)||_{L^\infty(Q)} +  \bar{C}_p \gamma_N \frac{\mu}{\tau} ||\Theta(\phi_1) - \Theta(\phi_2)||_{L^\infty(Q)}.
\end{eqnarray*}

Then, by Lemmas \ref{PropertiesEtcFirst} and \ref{icontLp01}, we finally have
\begin{eqnarray*}
  ||\Psi(l, \phi_1) - \Psi(l, \phi_2)||_{L^\infty(Q)} \leq C ||\phi_1 - \phi_2||_{L^\infty(Q)},
\end{eqnarray*}
where $C$ depends on $\bar{C}_p$, $C_1$, $C_2$, $\gamma$, $\gamma_N$, $\mu$, $\tau$ and the immersion constant.

To show that $\Psi (l, \cdot) $ is compact, we use the fact that the immersion
$W^{2,1}_4(Q) \hookrightarrow L^\infty(Q)$ is compact and that
$\Psi(l, \cdot)$ is the composition between the inclusion operator and the
solution operator, i.e., $\Psi(l, \cdot): L^\infty(Q)
\rightarrow W^{2,1}_4(Q) \rightarrow L^\infty(Q)$.
\hfill$\Box$
\\

\begin{lemma}
\label{unifcon1}
Given a bounded subset $B\subset L^\infty(Q)$, for each $\phi \in B$, the mapping $\Psi(\cdot, \phi): [0, 1] \rightarrow L^\infty(Q)$ is uniformly continuous  with respect to $B$.
\end{lemma}

\noindent{\bf Proof:}
Since $B \in L^\infty(Q)$ is bounded, there is $r_B \ge 0$ such that, for any $\phi \in B$, we have $||\phi||_{L^\infty(Q)} \leq r_B$. Now, let us fix $\phi \in L^\infty(Q)$ and consider $l_1, l_2 \in [0,1]$ and denote $\Psi(l_1, \phi) = D_1$, $\Psi(l_2, \phi) = D_2$ and $\tilde{D} = D_1 - D_2$. Then, $\tilde{D}$ satisfies
\begin{equation}
\label{B}
\left\{
\begin{array}{lcl}
\displaystyle \frac{\partial \tilde{D}}{\partial t} - \sigma \Delta \tilde{D} + \gamma l_2\Lambda(\phi)  \tilde{D} +  \gamma_N l_2 \Theta(\phi) \tilde{D} + \tau \tilde{D} = \\
\gamma\Lambda(\phi)D_1 (l_1 - l_2) - \gamma_N \Theta(\phi)D_1(l_1 - l_2), &\textup{in}& Q,
\vspace{0.2cm}
\\
\displaystyle \frac{\partial \tilde{D}}{\partial \eta} = 0, &\textup{on}& \Gamma,
\vspace{0.2cm}
\\
\displaystyle \tilde{D}(\cdot, 0) = \tilde{D}_0(\cdot) = 0, &\textup{in}& \Omega.
\end{array}
\right.
\end{equation}

Using the Proposition \ref{sol. Neumann} and the fact that $L^{\infty}(Q) \hookrightarrow L^4(Q)$, $D_1 \leq \frac{\mu}{\tau}$, we get
\begin{eqnarray*}
  ||\tilde{D}||_{W_{4}^{2,1}(Q)} \leq C_p ||\gamma\Lambda(\phi)D_1 (l_1 - l_2) - \gamma_N \Theta(\phi)D_1(l_1 - l_2)||_{L^4(Q)} \\
  \leq \bar{C}_p  \gamma \frac{\mu}{\tau} |l_1 - l_2| ||\Lambda(\phi)||_{L^\infty(Q)} + \bar{C}_p \gamma_N \frac{\mu}{\tau} |l_1 - l_2| ||\Theta(\phi)||_{L^\infty(Q)}.
\end{eqnarray*}

Then, by Lemmas \ref{PropertiesEtcFirst} and \ref{icontLp01}, we finally have
\begin{eqnarray*}
  ||\Psi(l_1,\phi) - \Psi(l_2,\phi)||_{L^\infty(Q)} \leq C |l_1 - l_2|,
\end{eqnarray*}
where $C$ depends on $\bar{C}_p$, $\gamma$, $\gamma_N$, $\mu$, $\tau$, $r_N$, $T$, $C_{\lambda}$, $||N_0||_{L^\infty(\Omega)}$, $||A_0||_{L^\infty(\Omega)}$ and the immersion constant.
\hfill$\Box$
\\

\begin{lemma}
\label{estimativa1}
Suppose $D_0 \leq \frac{\mu}{\tau}$ a.e. in $\Omega$, then there exists a number $\rho > 0$ such that, for any $l \in [0,1]$ and any possible
fixed point $D \in L^\infty(Q)$ of $\Psi(l, \cdot)$, there holds $\|D\|_{L^\infty(Q)} < \rho$.
\end{lemma}

\noindent {\bf Proof:} Let $D \in L^\infty(Q)$ such that $\Psi(l, D) = D$. The analogous demonstration made in Proposition \ref{nãonegativa1} guarantees us $||D||_{L^\infty(Q)} \leq \frac{\mu}{\tau}$. Therefore, just take $\rho = \frac{\mu}{\tau} + 1$.
\hfill$\Box$
\\

\begin{lemma}
\label{fix1}
The mapping $\Psi(0,\cdot): L^\infty(Q) \rightarrow L^\infty(Q)$ has a unique fixed point.
\end{lemma}

\noindent {\bf Proof:}
Indeed, letting $l =0$ in \ref{P6}, $D$ is a fixed point of $\Psi(0, \cdot)$ if, and only if,
$D$ is the unique solution to the problem
\begin{equation*}
\label{P0}
\left\{
\begin{array}{lcl}
\displaystyle	
\frac{\partial D}{\partial t} = \sigma \Delta D + \mu \chi_{\omega} - \tau D, &\textup{in}& Q,
\vspace{0.2cm}
\\
\displaystyle \frac{\partial D}{\partial \eta} (\cdot) =0, &\textup{on}& \Gamma,
\vspace{0.2cm}
\\
\displaystyle D(\cdot,0) = D_{0} (\cdot), &\textup{in}& \Omega.
\end{array}
\right.
\end{equation*}

But Proposition {\ref{sol. Neumann}} guarantees the  existence of a unique solution $D \in W^{2,1}_{4}(Q)\hookrightarrow L^{\infty}(Q)$
of this last problem; therefore $\Psi(0, \cdot )$ has a unique fixed point in $L^\infty(Q)$.
\hfill$\Box$
\\

\begin{proposition}
\label{existenciaA}
There is a nonnegative solution $\hat{D} \in W^{2,1}_4(Q)$ of the problem (\ref{P3}).
\end{proposition}

\noindent {\bf Proof:}
From Lemmas \ref{lemaAi}, \ref{con1}, \ref{unifcon1}, \ref{estimativa1} and \ref{fix1}, we conclude that the mapping $\Psi: [0,1] \times L^\infty(Q) \rightarrow L^\infty(Q)$ satisfies the hypotheses of the Leray-Schauder's fixed point theorem (see Friedman \cite[pp.~189, Theorem 3]{Friedman}). Thus, there exists $\hat{D} \in L^\infty(Q)$ such that $\Psi(1, \hat{D}) = \hat{D}$. Moreover, by Lemmas \ref{lemaAi} and \ref{nãonegativa1}, $\hat{D} \in W^{2,1}_4(Q)$ is nonnegative and $\hat{D}$ is the required solution of (\ref{P3}).
\hfill$\Box$
\\

\section{Proof of Theorem \ref{Teorema1}}

\begin{proposition}
\label{h2}
There is a nonnegative solution $(\hat{N}, \hat{A}, \hat{D}) \in L^\infty(Q) \times L^\infty(Q) \times W^{2,1}_4(Q)$ of the modified problem (\ref{P01}).
\end{proposition}

\noindent {\bf Proof:}
Just combine the Proposition \ref{existenciaA}, the Remark \ref{obs1} and the Lemma \ref{PropertiesEtcFirst}.
\hfill$\Box$
\\

\begin{remark}
\label{estimativaN}
We affirm that $\hat{N}, \hat{A} \in W$. Indeed, by Lemma \ref{PropertiesEtcFirst} we know that $\hat{N} = \Theta(\hat{D}), \hat{A} = \Lambda(\hat{D}) \in L^\infty(Q)$. Moreover, returning to the first equation of (\ref{P01}), using the Lemmas \ref{PropertiesEtcFirst} and \ref{nãonegativa1}, it follows that:
\begin{equation}\label{Nt}
\begin{array}{cccc}
\displaystyle
\bigg|\frac{\partial \hat{N}}{\partial t}\bigg| \leq r_N + \mu_N (||N_0||_{L^\infty(\Omega)} + r_N T) + \beta_1 (||N_0||_{L^\infty(\Omega)} + r_N T) C_{\lambda}||A_0||_{L^\infty(\Omega)} \\
\displaystyle
+ \alpha_N \gamma_N \frac{\mu}{\tau} (||N_0||_{L^\infty(\Omega)} + r_N T),
\end{array}
\end{equation}
a.e. in $Q$, i.e., $\hat{N}_t \in L^\infty(Q)$.

Moreover, returning to the second equation of (\ref{P01}) and using, again, the Lemmas \ref{PropertiesEtcFirst} and \ref{nãonegativa1}, we get:
\begin{equation}\label{At}
\begin{array}{cccc}
\displaystyle
\bigg|\frac{\partial \hat{A}}{\partial t}\bigg| \leq  r_A C_{\lambda}||A_0||_{L^\infty(\Omega)} + \frac{r_A}{k_A}(C_{\lambda}||A_0||_{L^\infty(\Omega)})^2
+ (\mu_A+\epsilon_A)C_{\lambda}||A_0||_{L^\infty(\Omega)} \\
\displaystyle
+ \alpha_A\gamma_A \frac{\mu}{\tau} C_{\lambda}||A_0||_{L^\infty(\Omega)},
\end{array}
\end{equation}
a.e. in $Q$, i.e., $\hat{A}_t \in L^\infty(Q)$.
\end{remark}

\begin{proposition}
\label{j1}
There is a nonnegative solution $(N, A, D) \in W \times W \times W^{2,1}_4(Q)$ of problem (\ref{0riginalEquations}).
\end{proposition}

\noindent {\bf Proof:}
Just combine the Proposition \ref{h2} and the Remarks \ref{obs2} and \ref{estimativaN}.
\hfill$\Box$
\\

\begin{proposition}
\label{j2}
The solution $(N, A, D)$ of the problem (\ref{0riginalEquations}) is unique.
\end{proposition}

\noindent {\bf Proof:}
Let $(N_1, A_1, D_1)$ and $(N_2, A_2, D_2)$ be solutions to the problem (\ref{0riginalEquations}); if $\tilde{N} = N_1 - N_2, \tilde{A} = A_1 - A_2$ and $\tilde{D} = D_1 - D_2$, then $\tilde{N}$, $\tilde{A}$ and $\tilde{H}$ satisfy the following problems, respectively:
\begin{equation}
\label{original1}
\left\{
\begin{array}{lcl}
\displaystyle
\frac{\partial \tilde{N}}{\partial t} = - \mu_N \tilde{N}  - \beta_1 A_1 \tilde{N} -\beta_1 N_2\tilde{A} - \alpha_N\gamma_N N_1 \tilde{D}-\alpha_N\gamma_N  D_2\tilde{N}, & \textup{in}& Q,
\vspace{0.2cm}
\\
\displaystyle
\displaystyle \tilde{N}(\cdot,0) = \tilde{N}_0(\cdot) = 0, &\textup{in}& \Omega,
\end{array}
\right.
\end{equation}

\begin{equation}
\label{original2}
\left\{
\begin{array}{lcl}
\displaystyle
\frac{\partial \tilde{A}}{\partial t} = r_A \tilde{A} -\frac{r_A}{k_A}(A_1 + A_2)\tilde{A}-(\mu_A+\epsilon_A)\tilde{A} - \alpha_A\gamma_A A_1 \tilde{D}- \alpha_A\gamma_A D_2\tilde{A}, &\textup{in}& Q,
\vspace{0.2cm}
\\
\displaystyle \tilde{A}(\cdot,0) = \tilde{A}_0(\cdot) = 0, &\textup{in}& \Omega,
\end{array}
\right.
\end{equation}

\begin{equation}
\label{original3}
\left\{
\begin{array}{lcl}
\displaystyle
\frac{\partial \tilde{D}}{\partial t} = \sigma \Delta \tilde{D} - \gamma A_1 \tilde{D} -\gamma D_2\tilde{A} - \gamma_N N_1 \tilde{D} -\gamma_N D_2\tilde{N} - \tau \tilde{D}, &\textup{in}& Q,
\vspace{0.2cm}
\\
\displaystyle \frac{\partial \tilde{D}}{\partial \eta} (\cdot) =0, &\textup{on}& \Gamma,
\vspace{0.2cm}
\\
\displaystyle \tilde{D}(\cdot,0) = \tilde{D}_0(\cdot) = 0, &\textup{in}& \Omega.
\end{array}
\right.
\end{equation}

Multiplying the first equation of (\ref{original1}) by $\tilde{N}$, integrating into $\Omega$, using the fact that $N_1 \leq ||N_0||_{L^\infty(\Omega)} + r_N T$ and the inequality of Young, we have
\begin{eqnarray*}
\frac{1}{2}\frac{d}{dt} \int_{\Omega} \tilde{N}^2 dx &=& -\mu_N \int_{\Omega} \tilde{N}^2 dx -\beta_1 \int_{\Omega} N_2\tilde{A}\tilde{N} dx - \alpha_N\gamma_N \int_{\Omega} N_1\tilde{D}\tilde{N} dx \\
\\
&-& \alpha_N\gamma_N \int_{\Omega} D_2 \tilde{N}^2 dx \\
\\
&\leq& (||N_0||_{L^\infty(\Omega)} + r_N T) \bigg(\beta_1 \int_{\Omega}|\tilde{A}||\tilde{N}| dx + \alpha_N\gamma_N \int_{\Omega} |\tilde{D}||\tilde{N}| dx\bigg) \\
&\leq& C \int_{\Omega} (\tilde{A}^2 + \tilde{N}^2 + \tilde{H}^2) dx,
\end{eqnarray*}
where $C$ depends on $\beta_1$, $\alpha_N$, $\gamma_N$, $r_N$, $T$ and $||N_0||_{L^\infty(\Omega)}$.

Now, multiplying the first equation of (\ref{original2}) by $\tilde{A}$, integrating into $\Omega$, using the fact that $A_1 \leq C_{\lambda}||A_0||_{L^\infty(\Omega)} $ and the inequality of Young, we obtain
\begin{eqnarray*}
  \frac{1}{2}\frac{d}{dt} \int_{\Omega} \tilde{A}^2 dx &=& r_A \int_{\Omega}\tilde{A}^2 dx - \frac{r_A}{k_A}\int_{\Omega}(A_1 + A_2)\tilde{A}^2 dx -(\mu_A+\epsilon_A)\int_{\Omega}\tilde{A}^2 dx \\
&-& \alpha_A \gamma_A \int_{\Omega} A_1 \tilde{D}\tilde{A} dx - \alpha_A \gamma_A \int_{\Omega}D_2\tilde{A}^2 dx \\
&\leq& r_A \int_{\Omega}|\tilde{A}|^2 dx + \alpha_A \gamma_A C_{\lambda}||A_0||_{L^\infty(\Omega)} \int_{\Omega} |\tilde{D}||\tilde{A}| dx \\
&\leq& C \int_{\Omega} (\tilde{A}^2 + \tilde{N}^2 + \tilde{H}^2) dx,
\end{eqnarray*}
where $C$ depends on $r_A$, $\alpha_A$, $\gamma_A$ $C_{\lambda}$ and $||A_0||_{L^\infty(\Omega)}$.

Lastly, multiplying the first equation of (\ref{original3}) by $\tilde{H}$, integrating into $\Omega$, using the fact that $D \leq \frac{\mu}{\tau}$ and the inequality of Young, we obtain
\begin{eqnarray*}
  \frac{1}{2}\frac{d}{dt} \int_{\Omega} \tilde{D}^2 dx &=& -\sigma \int_{\Omega} |\nabla \tilde{D}|^2 dx - \gamma \int_{\Omega} A_1 {\tilde{D}}^2 dx - \gamma \int_{\Omega} D_2\tilde{A} \tilde{D} dx \\
&-& \gamma_N \int_{\Omega}N_1 \tilde{D}^2 dx - \gamma_N \int_{\Omega} D_2 \tilde{N}\tilde{D} dx - \tau \int_{\Omega}{\tilde{D}}^2 dx \\
&\leq& \gamma\frac{\mu}{\tau} \int_{\Omega} |\tilde{A}| |\tilde{D}| dx + \gamma_N \frac{\mu}{\tau} \int_{\Omega}|\tilde{N}||\tilde{D}| dx \\
&\leq& C \int_{\Omega} (\tilde{A}^2 + \tilde{N}^2 + \tilde{H}^2) dx,
\end{eqnarray*}
where $C$ depends on $\gamma$, $\gamma_N$, $\mu$ and $\tau$.

Thus,
\begin{eqnarray*}
  \frac{d}{dt}\bigg( \int_{\Omega} (|\tilde{N}|^2 + |\tilde{A}|^2 + |\tilde{D}|^2) dx \bigg) &\leq& C \int_{\Omega} (|\tilde{N}|^2 + |\tilde{A}|^2+ |\tilde{D}|^2) dx,
\end{eqnarray*}
and using the Gronwall's inequality, we finally
\begin{equation*}
\int_{\Omega} (|\tilde{N}|^2 + |\tilde{A}|^2 + |\tilde{D}|^2) dx \leq e^{CT}  \int_{\Omega} (|\tilde{N}_0|^2 + |\tilde{A}_0|^2 + |\tilde{D}_0|^2) dx = 0,
\end{equation*}
that is, $||\tilde{N}(\cdot, t)||_{L^2(\Omega)}^2 + ||\tilde{A}(\cdot, t)||_{L^2(\Omega)}^2 + ||\tilde{D}(\cdot, t)||_{L^2(\Omega)}^2 = 0$, for all $t \in (0, T)$. Where we conclude $\tilde{N} = \tilde{A} = \tilde{D} = 0$ a.e. in $Q$ and therefore $N_1 = N_2, A_1 = A_2$ and $D_1 = D_2$ a.e. in $Q$.
\hfill$\Box$
\\

\section{Numerical simulations}

In this section, we provide numerical simulations illustrating different model behaviors. The settings and methods used to implement the simulations are the following. We consider the spatial domain as a square $\Omega=[0,L] \times [0,L]$, with $L=1$, discretized with $n=50$ steps $\Delta x = \Delta y = L/n=0.02$. The Laplacian $\Delta D$ is approximated by second order centered finite differences and the coupled ODE system arising from such discretization is solved with the method of lines in the software \textit{Mathematica}. The simulations run from time $t=0$ until $t=25$ (which is enough to achieve stationary behavior in all simulations).

The initial conditions for numerical simulations are $N(x,0)=N_2$, $A(x,0)=A_2$, $D(x,0)=0$, where $(N_2,A_2,0)$ is a globally asymptotically stable equilibrium point for the ODE system \eqref{início1} without treatment ($\nu=0$). The expressions for $N_2$ and $A_2$ are:
\[
N_2=\dfrac{r_N}{\mu_N+\beta_1 A_2},
\ \ \
A_2=\dfrac{r_A-\mu_A-\epsilon_A}{r_A}K_A.
\]
Such equilibrium is allways globally asymptotically stable in system \eqref{início1} (see details in \cite{Fassoni}). From the biological point of view, these initial conditions correspond to the start of chemotherapy application when a tumor is already a formed, where the normal cells were not able to control tumor growth, and no chemotherapy was applied until the tumor reached a stationary state.

To avoid large numbers and numerical instabilities, we re-scale the populations with respect to their possible maximum values, setting $N \leftarrow N/(r_N/\mu_N)$ and $A\leftarrow A/K_A$. Therefore, the population sizes range from $0$ to $1$. The re-scaled parameter values used in the model simulations were fixed to
\[
r_N=1, \
\mu_N=1, \
r_A = 1, \
K_A = 1, \
\beta_1 = 1.5, \
\mu_A =0.05, \
\epsilon_A =0.05, \
\]
\[
\tau_H=0.9, \
\gamma_N=0.1, \
\alpha_N =1, \
\gamma_A=1.
\]
These values were chosen to describe: normal cells that reach the equilibrium $N=r_N/\mu_N=1$ at absence of tumor cells; a tumor with the same carrying capacity of normal cells ($K_A=r_N/\mu_N=1$) and a greater absorption of the chemotherapeutic drug by tumor cells in comparison with normal cells ($\gamma_A>\gamma_N$), due to the drug specificity.

In order to illustrate different biological outcomes in the model simulations, we allowed the following parameters to assume different values: the chemotherapeutic drug cytotoxicity against cancer cells $\alpha_A$, the diffusion coefficient of the chemotherapeutic drug $\sigma$ and the chemotherapy infusion rate $\mu$. We will show that these properties of the drug and the infusion rate are crucial for determining an effective treatment. We also simulated different positions for the subset $\omega$, which is a mathematical description of a blood vessel crossing the tissue, from where the chemotherapy enters the tissue. The values for parameters $\alpha_A$, $\sigma$, $\mu$ and the position of $\omega$ used in each simulation are indicated in Table \ref{tableSims}. We present the following results.

\begin{center}
\begin{table}
\begin{tabular}{|c|c|c|c|c|c|c|}
\hline
Simulation & Figure & Outcome  & $\alpha_A$ & $\mu$ & $\sigma$ & $\omega$ \\
\hline
1  & \ref{fig1} & tumor persistence 			& 5 & 3 & 0.1 & $[0.45,0.55] \times [0.45,0.55]$ \\
\hline
2 & \ref{fig2} & tumor persistence 			& 10 & 3 & 0.1 & $[0,0.1]$ \\
\hline
3 & \ref{fig3} & tumor extinction 			& 10 & 6 & 0.1 & $[0,0.1]$ \\
\hline
4 & \ref{fig4} & tumor extinction 			& 10 & 3 & 0.2 & $[0,0.1]$ \\
\hline
5 & \ref{fig5} & tumor extinction 			& 20 & 3 & 0.1 & $[0,0.1]$ \\
\hline
\end{tabular}
\caption{Set-up of different simulations an their biological outcomes. Each row indicates the numerical values used for the chemotherapeutic parameters $\alpha_A$ (cytotoxicity), $\sigma$ (diffusion coefficient), $\mu$ (infusion rate), and the position of $\omega \subset \Omega \subset \mathbb{R}^2$. Simulation 1 was performed in a two-dimensional domain $\Omega=[0,1]\times [0,1]$, while simulations 2-5 were performed in a one-dimensional domain $\Omega=[0,1]$.}
\label{tableSims}
\end{table}
\end{center}

\begin{figure}[!htb]
\centering
\includegraphics[width=0.32\linewidth]{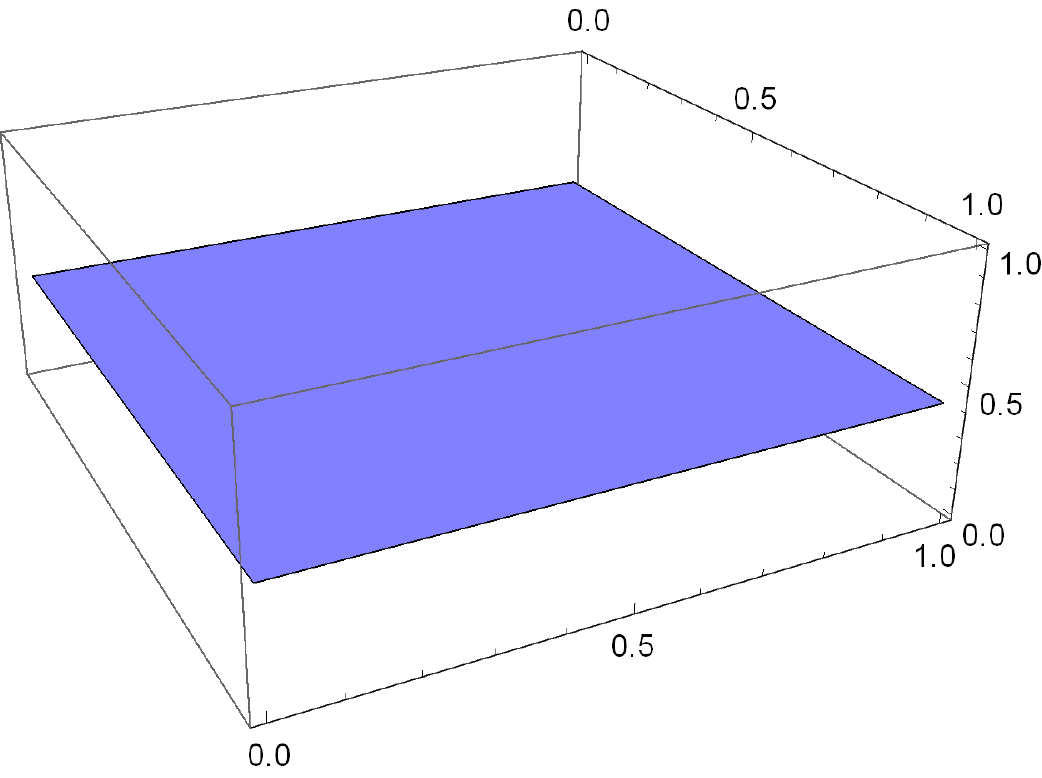}
\includegraphics[width=0.32\linewidth]{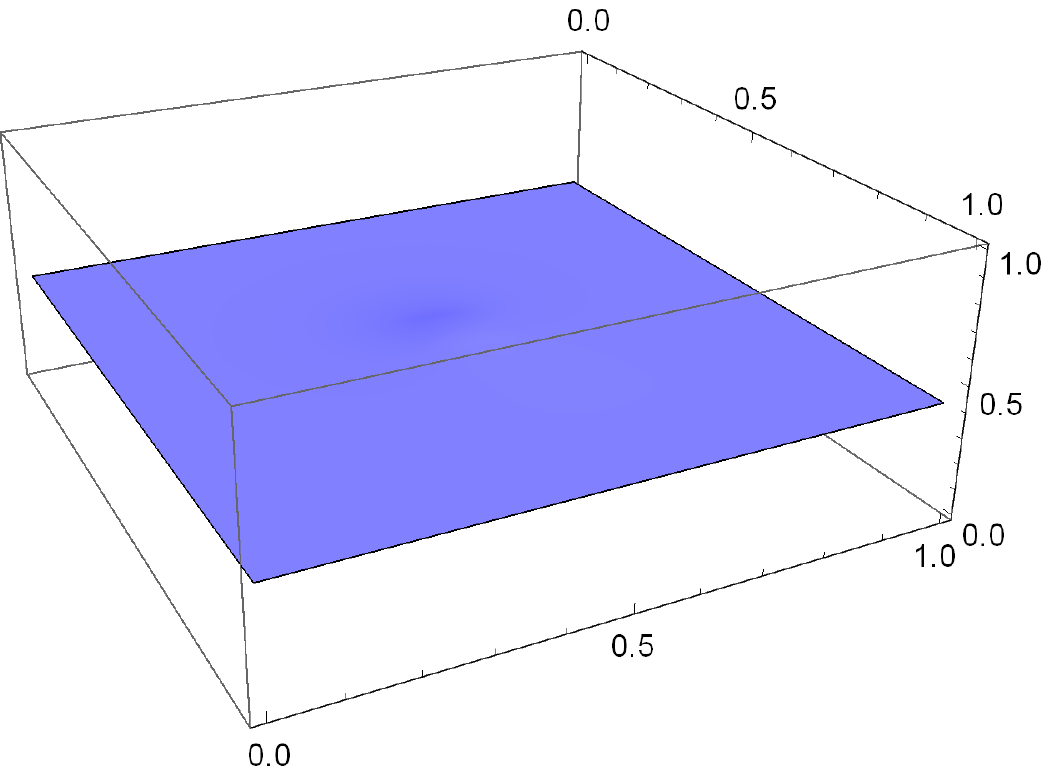}
\includegraphics[width=0.32\linewidth]{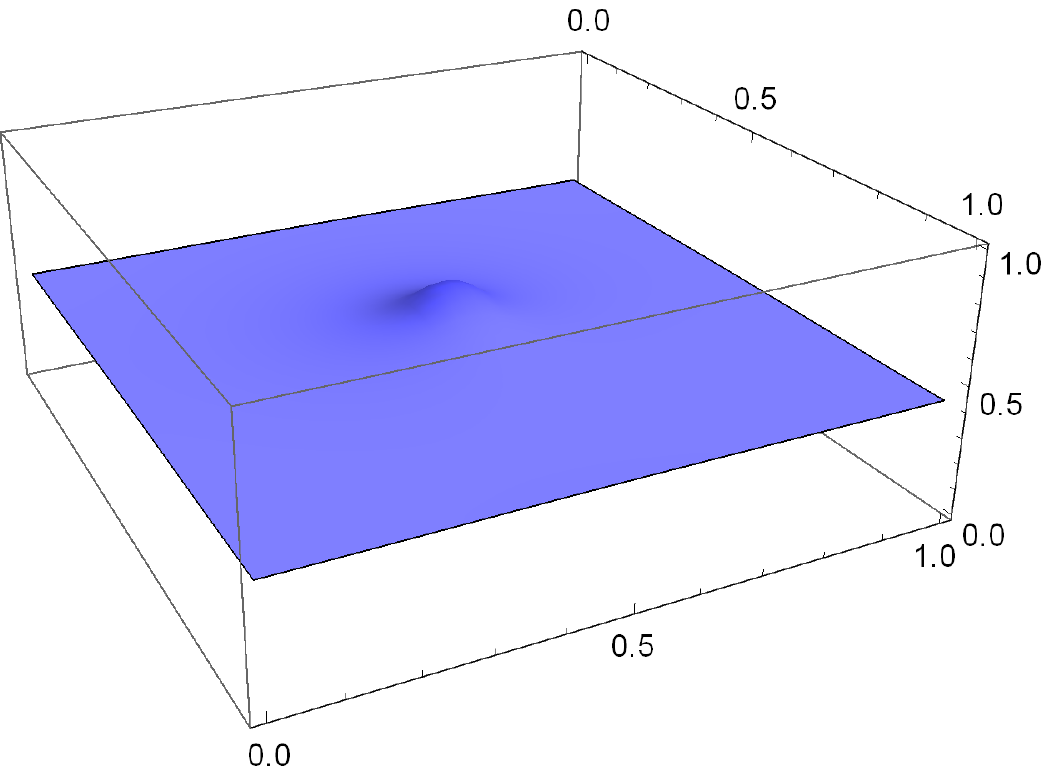}
\includegraphics[width=0.32\linewidth]{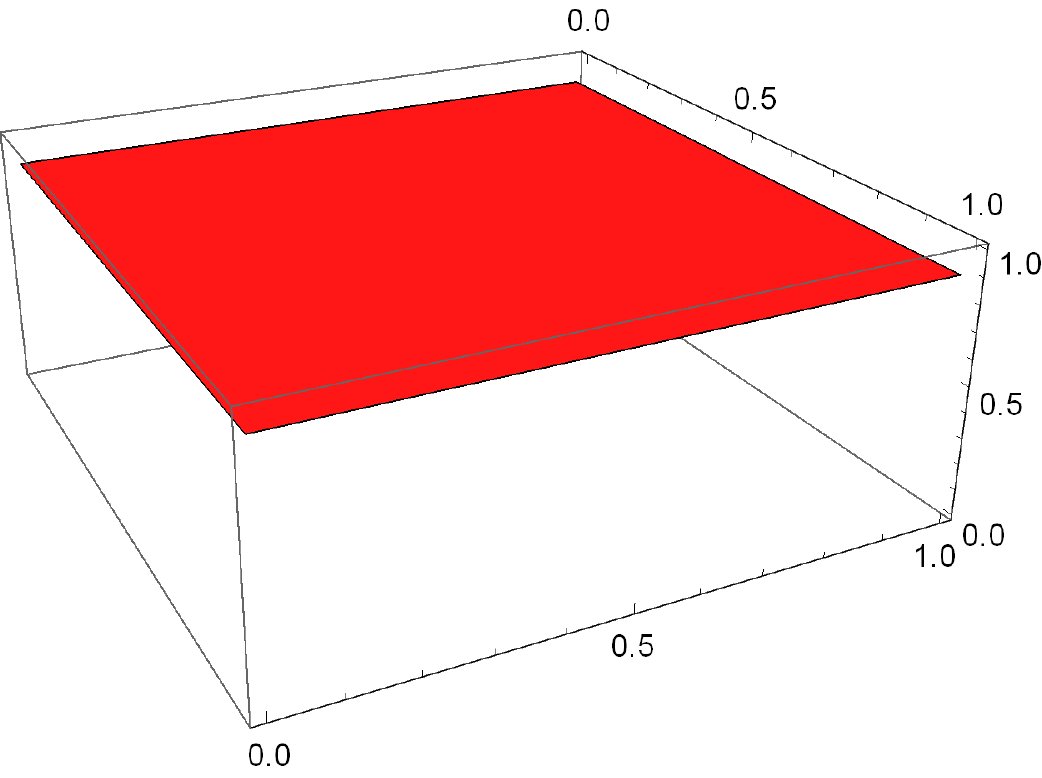}
\includegraphics[width=0.32\linewidth]{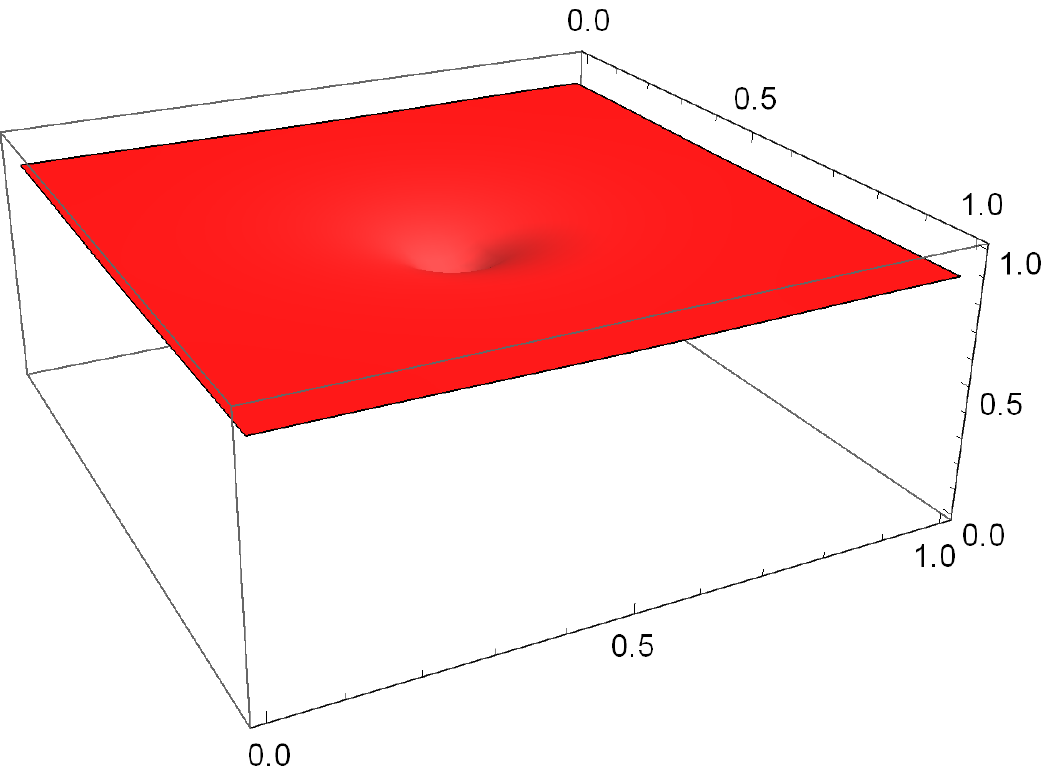}
\includegraphics[width=0.32\linewidth]{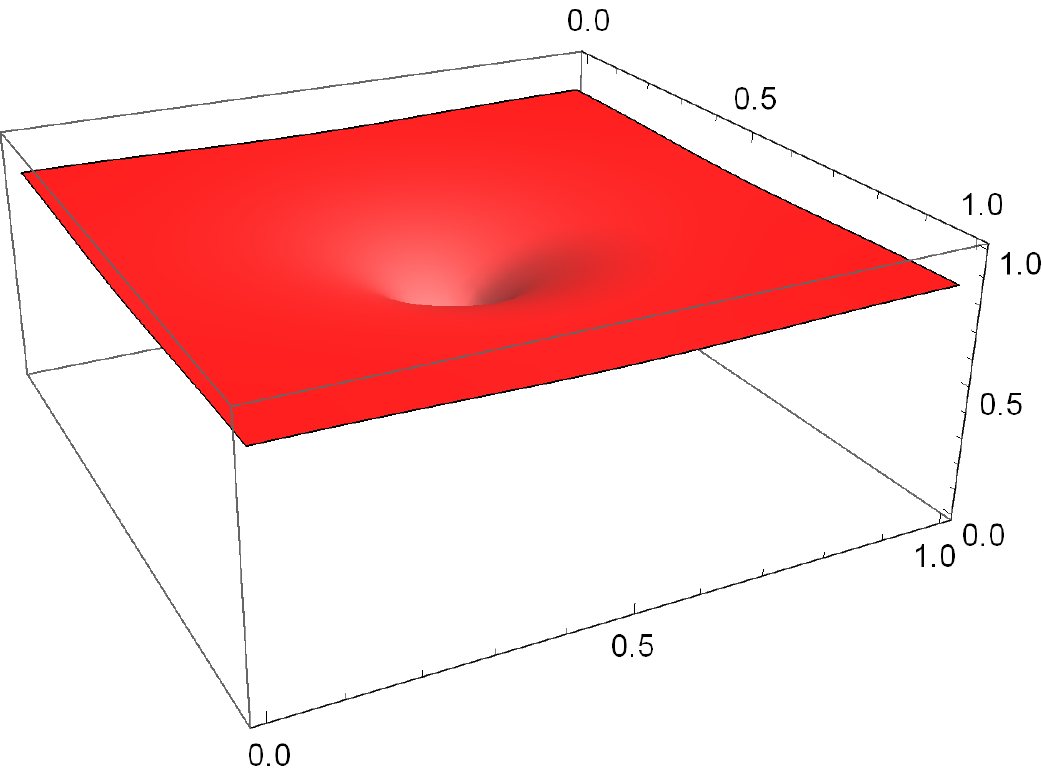}
\includegraphics[width=0.32\linewidth]{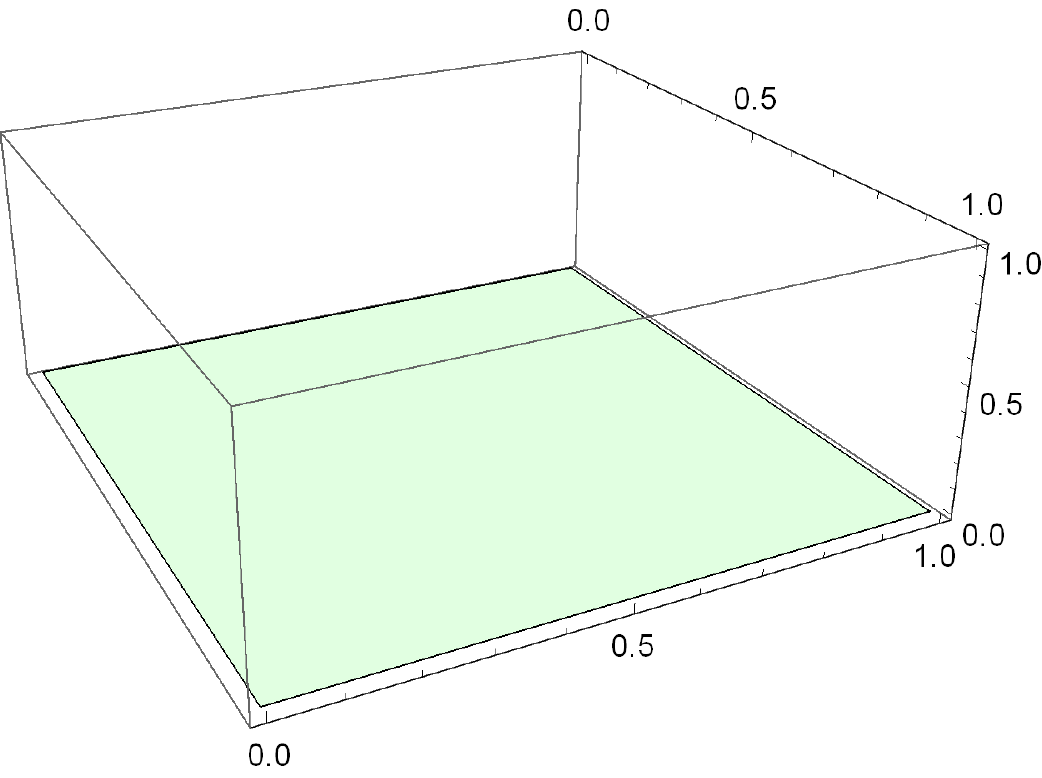}
\includegraphics[width=0.32\linewidth]{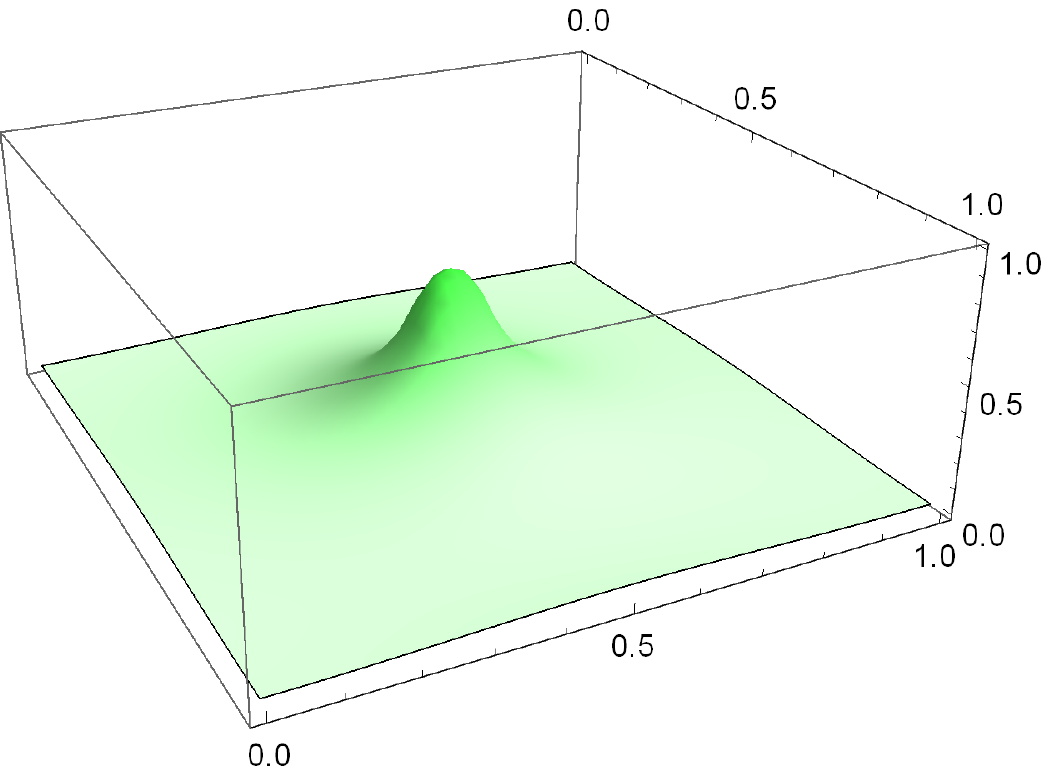}
\includegraphics[width=0.32\linewidth]{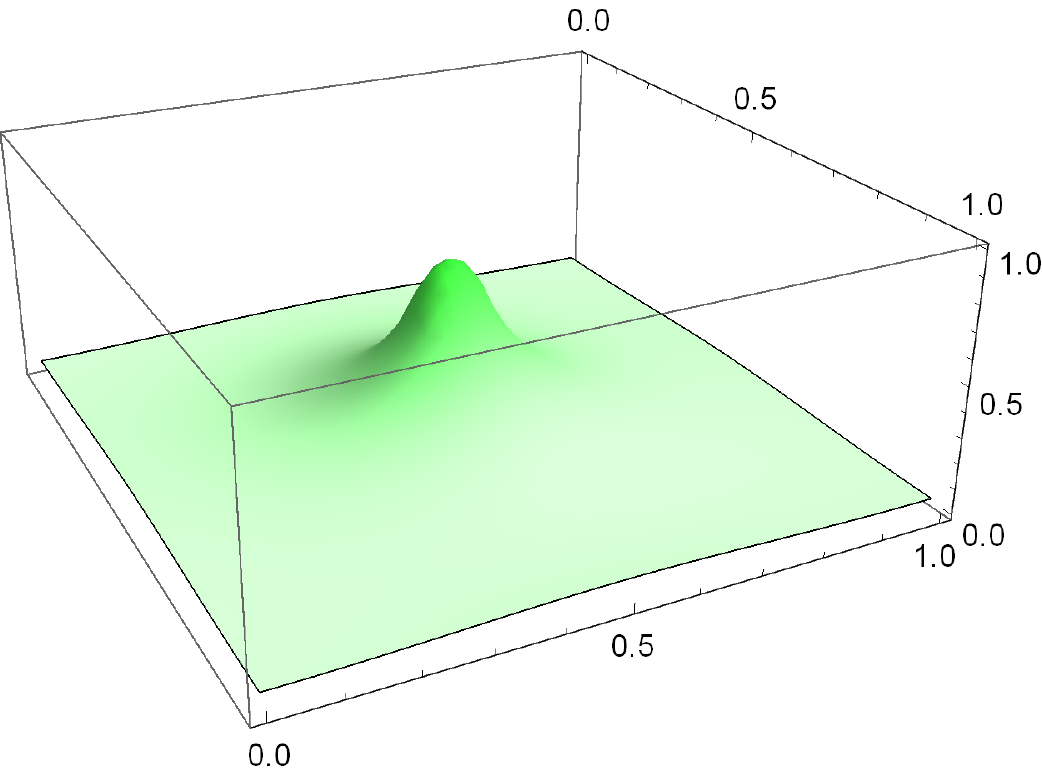}
\caption{Results of Simulation 1, for model (\ref{0riginalEquations}) within a two-dimensional domain $\Omega = [0,L] \times [0,L] =[0,1] \times [0,1]$. Plots of model solutions $N(x,y,t)$ (normal cells, blue, top row), $A(x,y,t)$ (cancer cells, red, middle row) and $D(x,y,t)$ (chemotherapeutic drug concentration, green, bottom row) at time points $t=0,1,15$ (columns 1,2 and 3, respectively). See Table \ref{tableSims} for parameter values used here. At time $t=0$, the tumor is spread trough the tissue, and as chemotherapy is applied ($t>0$), the tumor cells are reduced in the vicinity of the blood vessel, while the distant tumor cells  persist along time (the shape of the solution at time $t=15$ is stationary). Within the vicinity of the blood vessel, the removal of tumor cells allows the normal tissue to recover and grow.}
\label{fig1}
\end{figure}

In the first simulation of system \eqref{0riginalEquations}, we confirm that our model and numerical methods are able to reproduce the expected biological behavior (Figure \ref{fig1}). The blood vessel crosses the tissue at its center, i.e., $\omega =[0.45,0.55] \times [0.45,0.55]$. We use the following parameter values: $\alpha_A=5$, $\mu=3$, and $\sigma=0.1$. With such values, the chemotherapy is not able to lead to tumor extinction. We observe that tumor cells that are near the blood vessel are eliminated but not extinct by the chemotherapeutic effect, and those which are distant from the blood vessel persist (Figure \ref{fig1}).

\begin{figure}[!htb]
\centering
\includegraphics[width=0.32\linewidth]{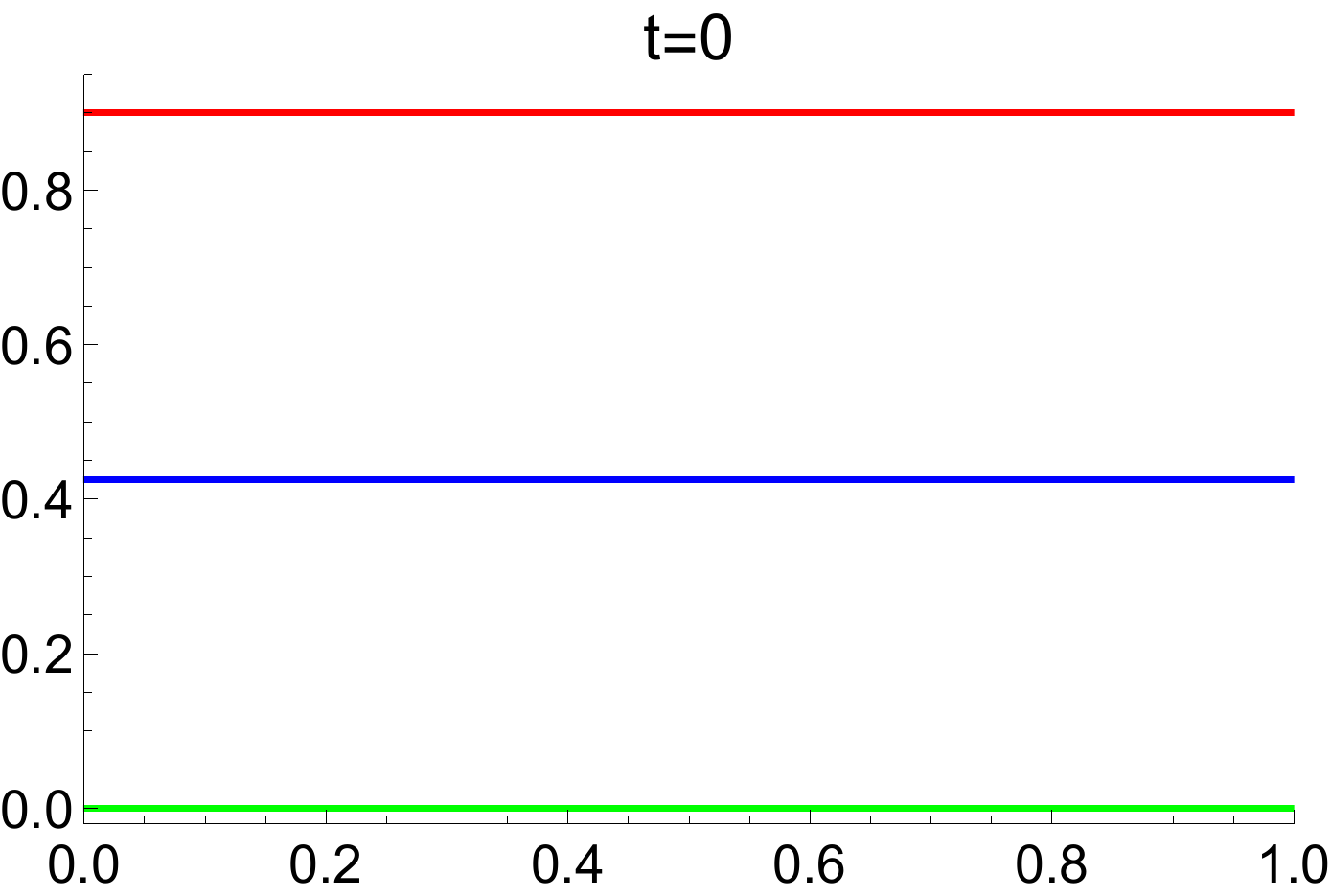}
\includegraphics[width=0.32\linewidth]{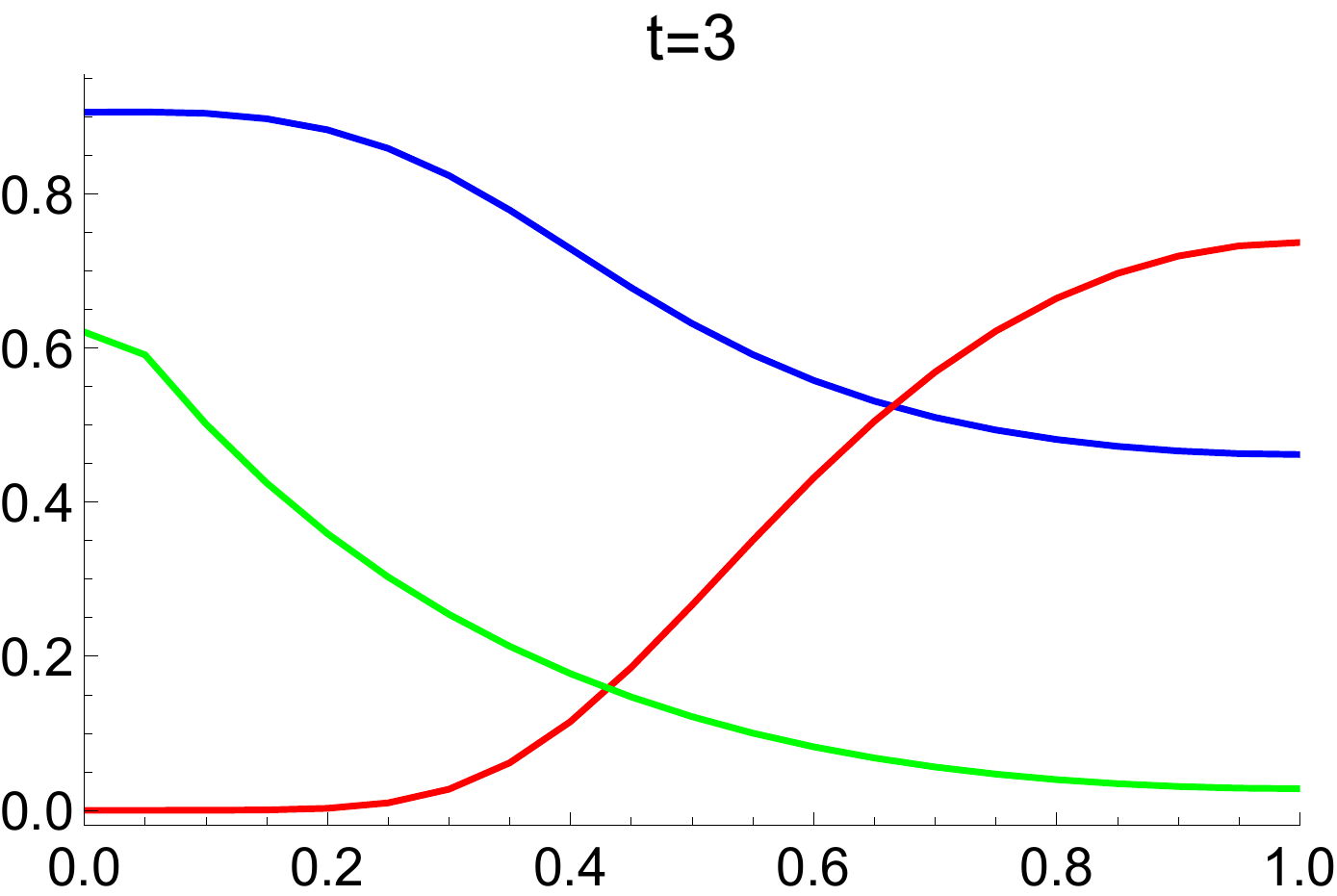}
\includegraphics[width=0.32\linewidth]{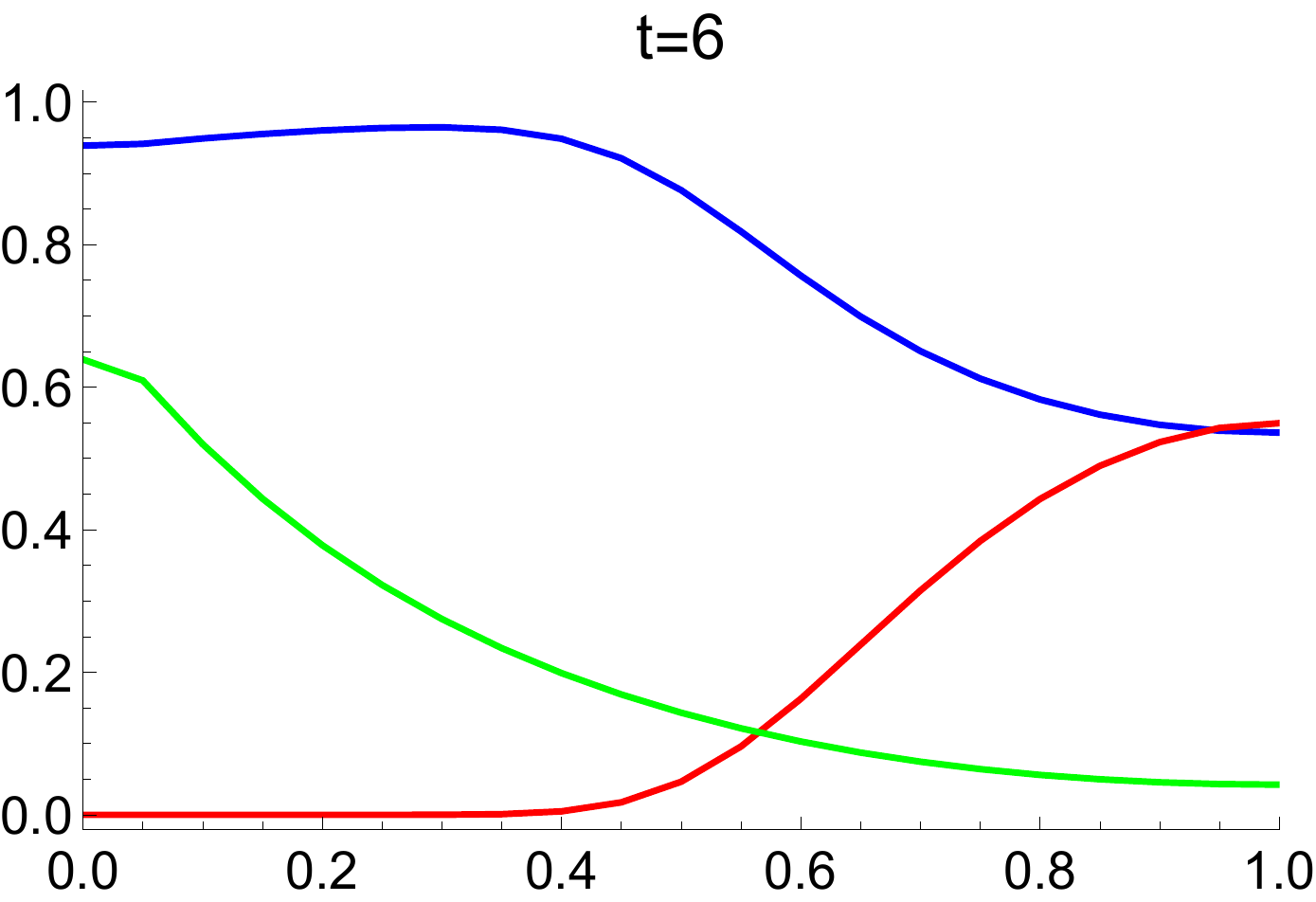}
\includegraphics[width=0.32\linewidth]{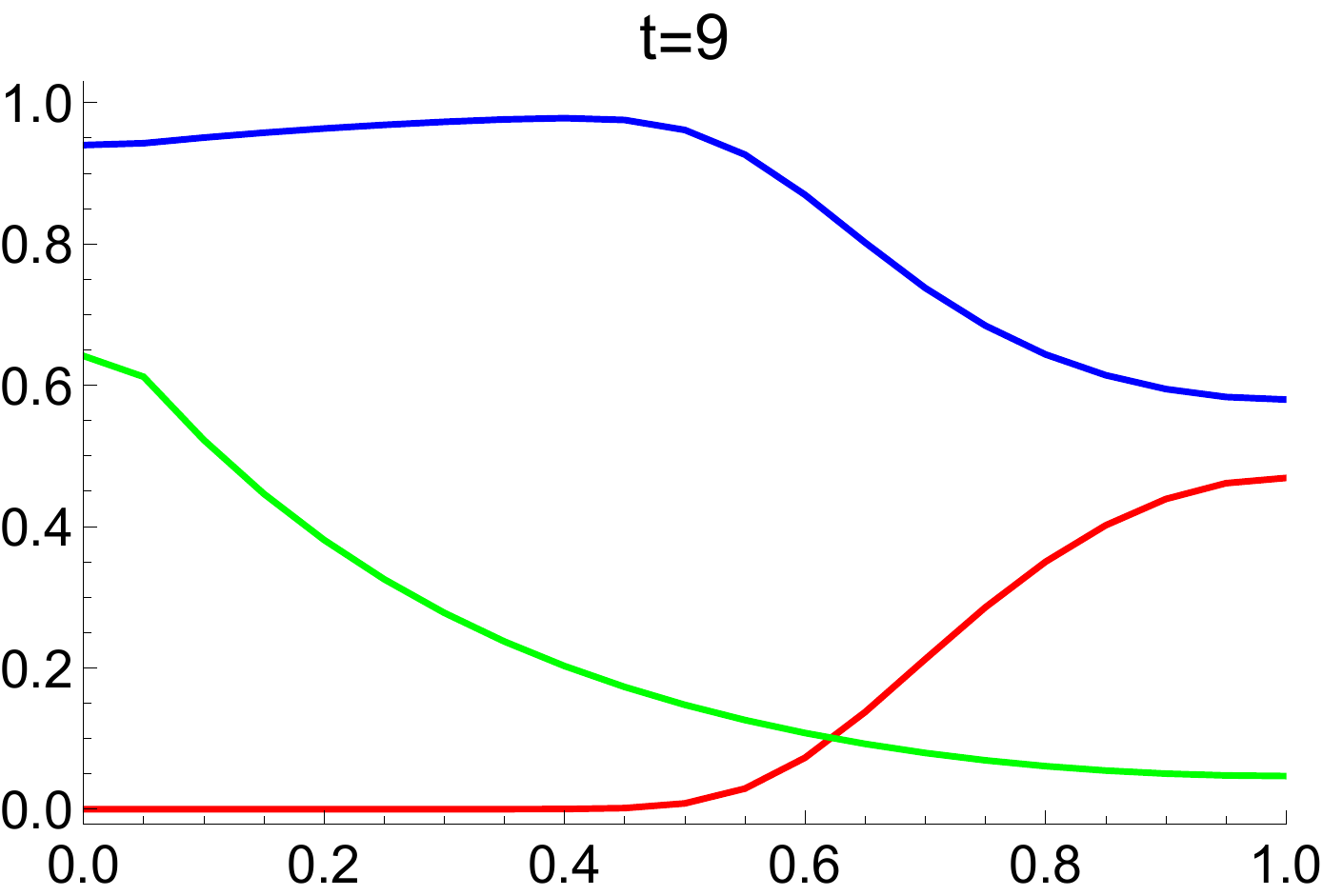}
\includegraphics[width=0.32\linewidth]{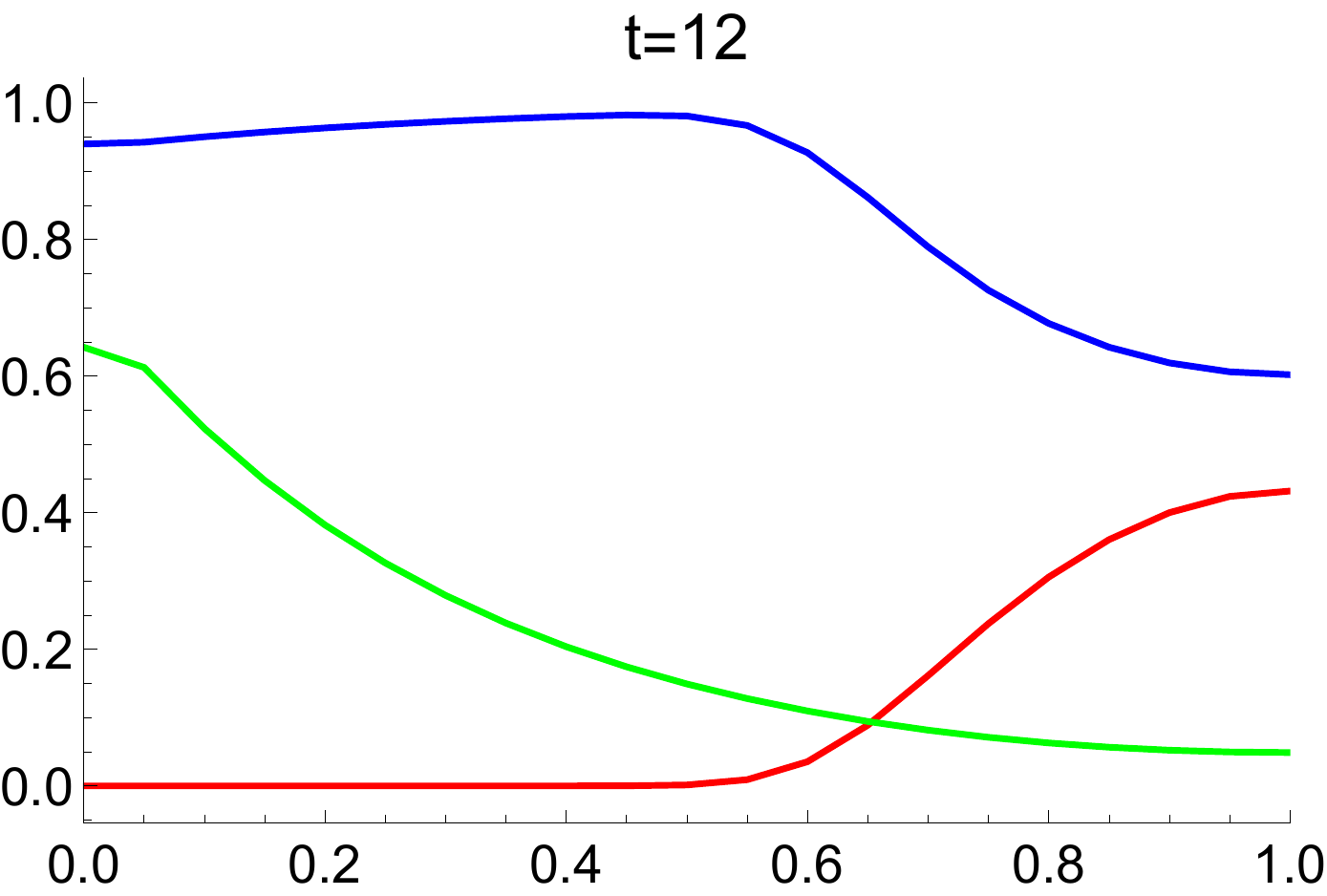}
\includegraphics[width=0.32\linewidth]{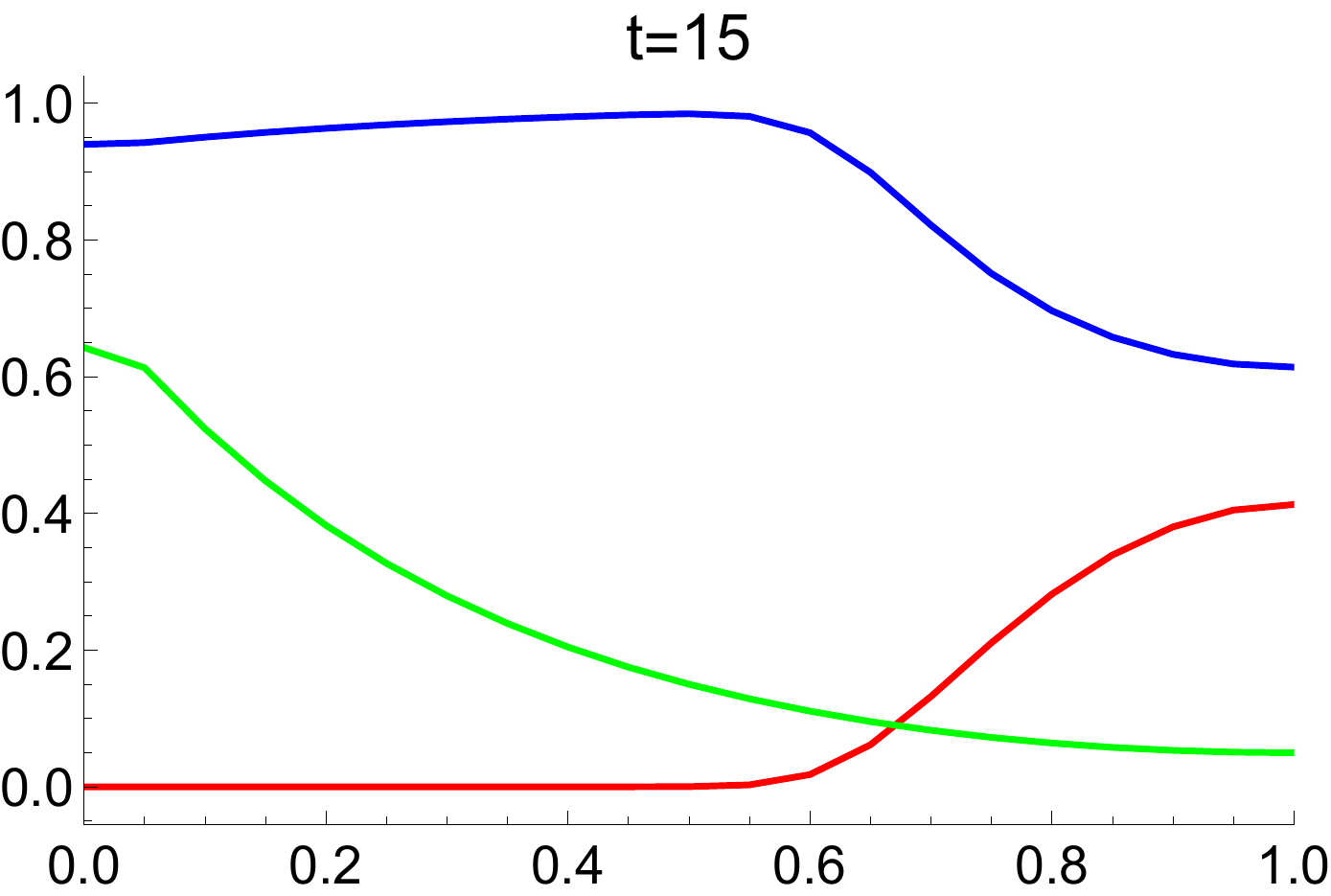}
\caption{Results of Simulation 2, with $\Omega = [0,L]= [0,1]$. Plots of model solutions $A(x,t)$ (cancer cells, red), $N(x,t)$ (normal cells, blue) and $D(x,t)$ (chemotherapeutic drug concentration, green) at time points $t=0,3,6,9,12,15$. See Table \ref{tableSims} for parameter values used here. At time $t=0$, the tumor is spread trough the tissue, and as chemotherapy is applied ($t>0$), the tumor cells are reduced and extinct within a given distance from the blood vessel ($x<0.6$), but not in the entire tissue ($x>0.6$). Within the region of tumor extinction, the removal of tumor cells release the normal tissue to recover and grow.}
\label{fig2}
\end{figure}

\begin{figure}[!htb]
\centering
\includegraphics[width=0.32\linewidth]{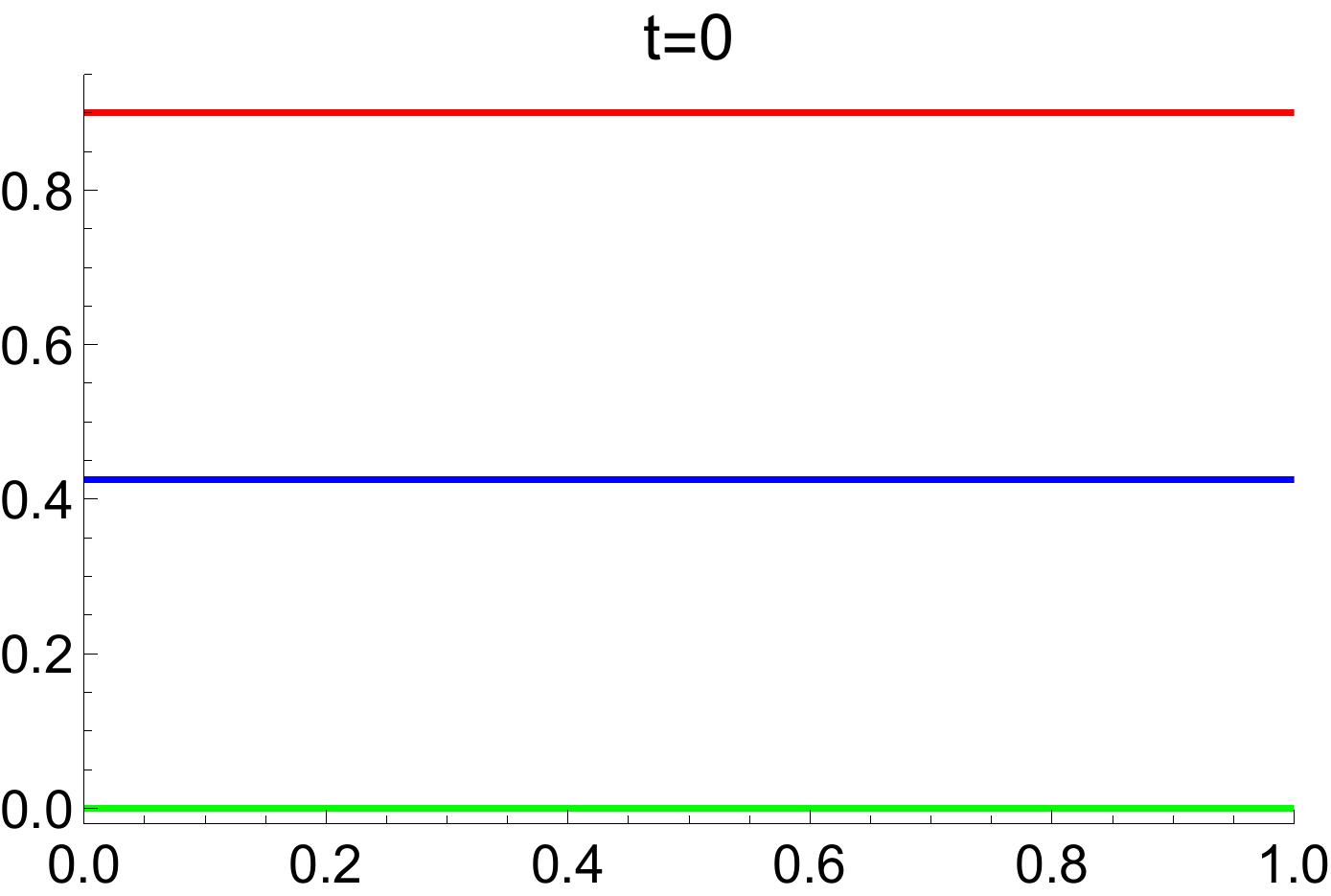}
\includegraphics[width=0.32\linewidth]{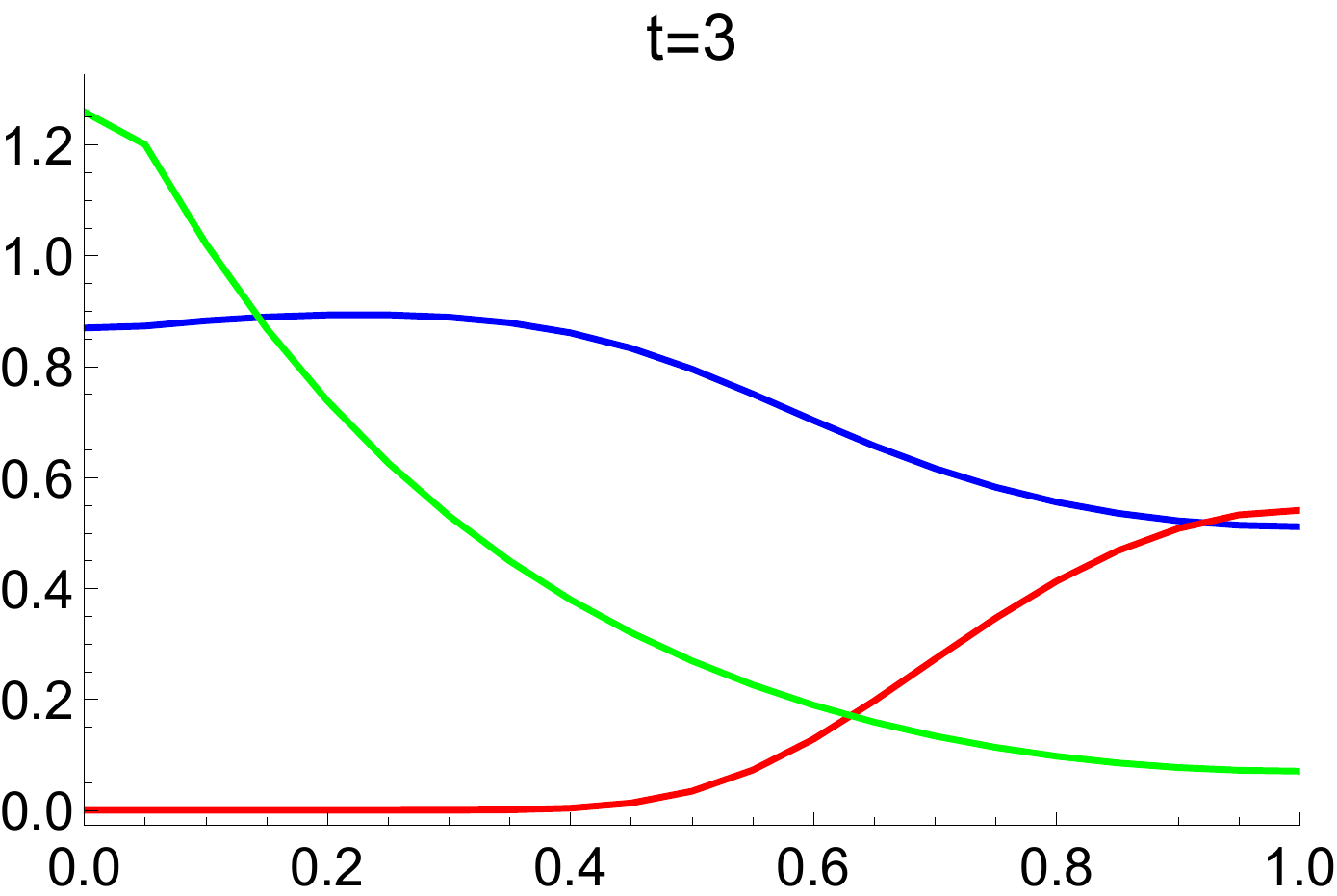}
\includegraphics[width=0.32\linewidth]{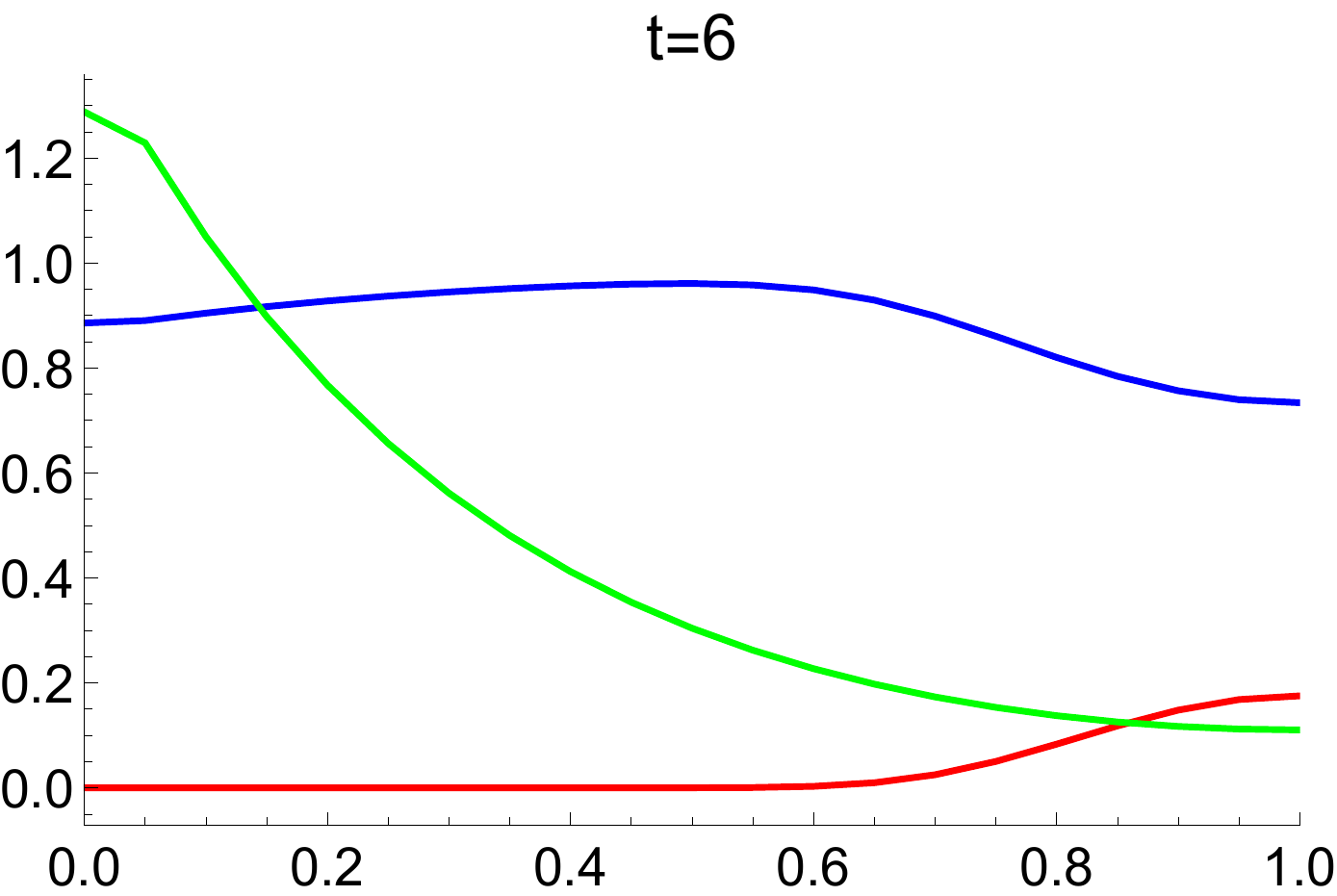}
\includegraphics[width=0.32\linewidth]{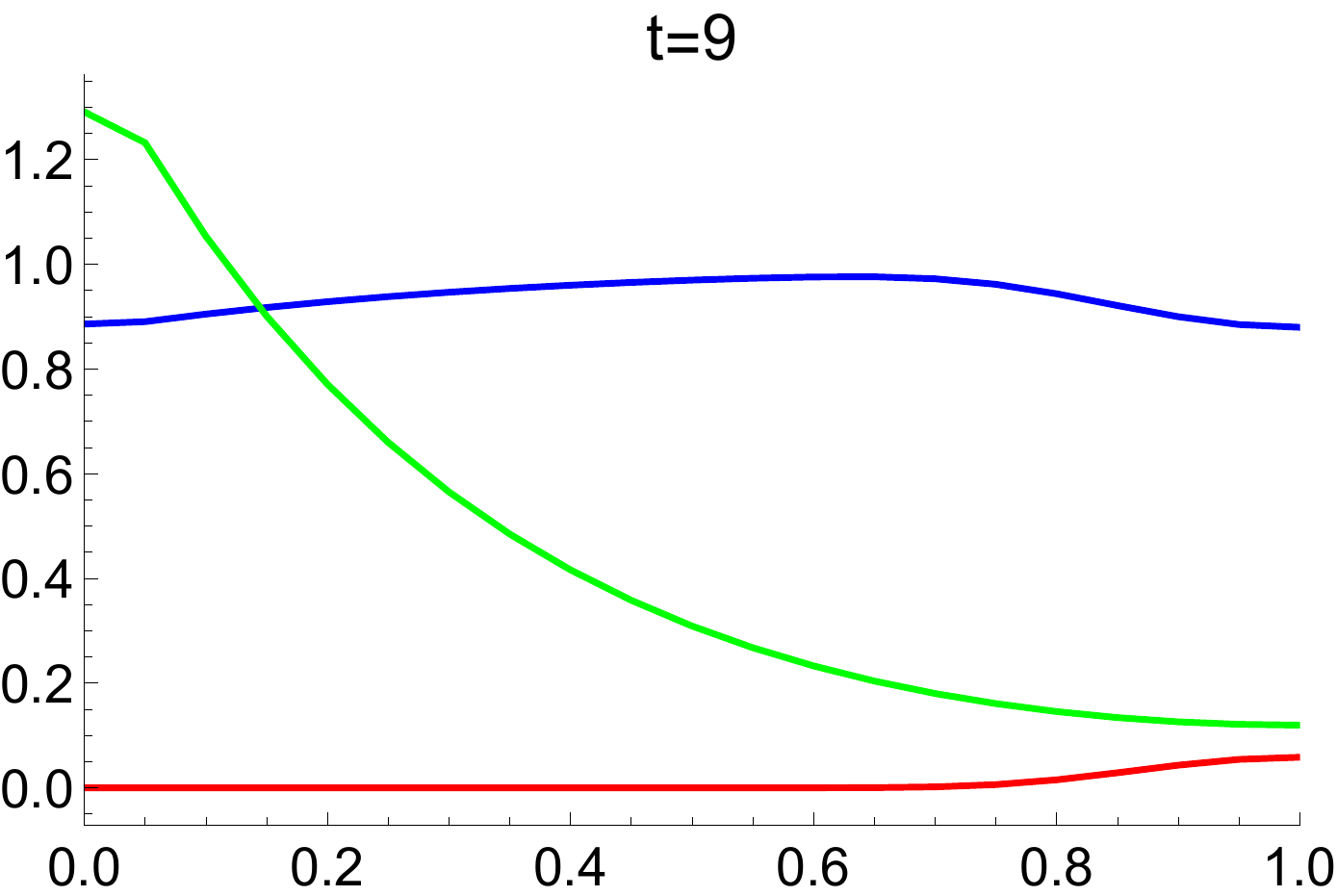}
\includegraphics[width=0.32\linewidth]{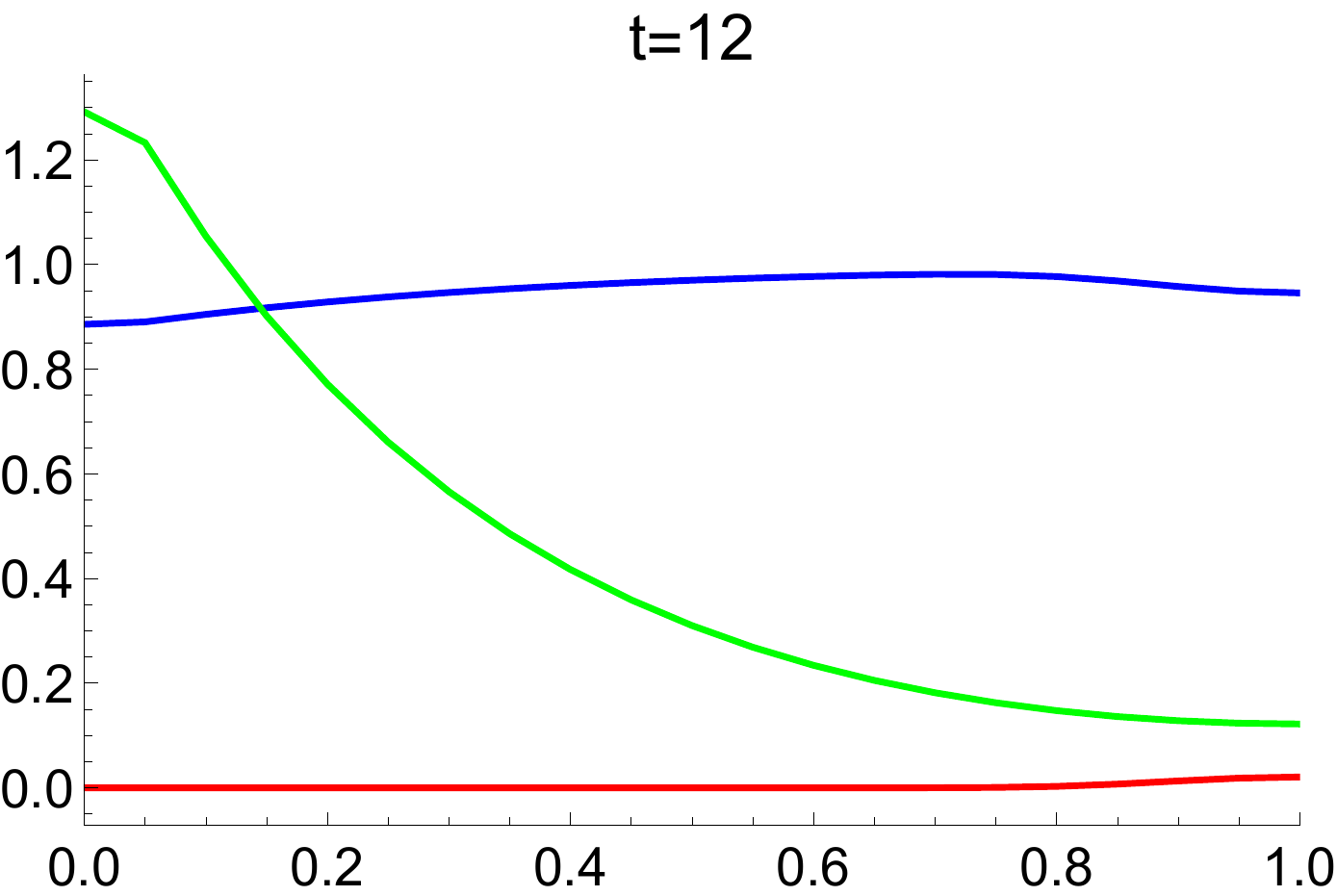}
\includegraphics[width=0.32\linewidth]{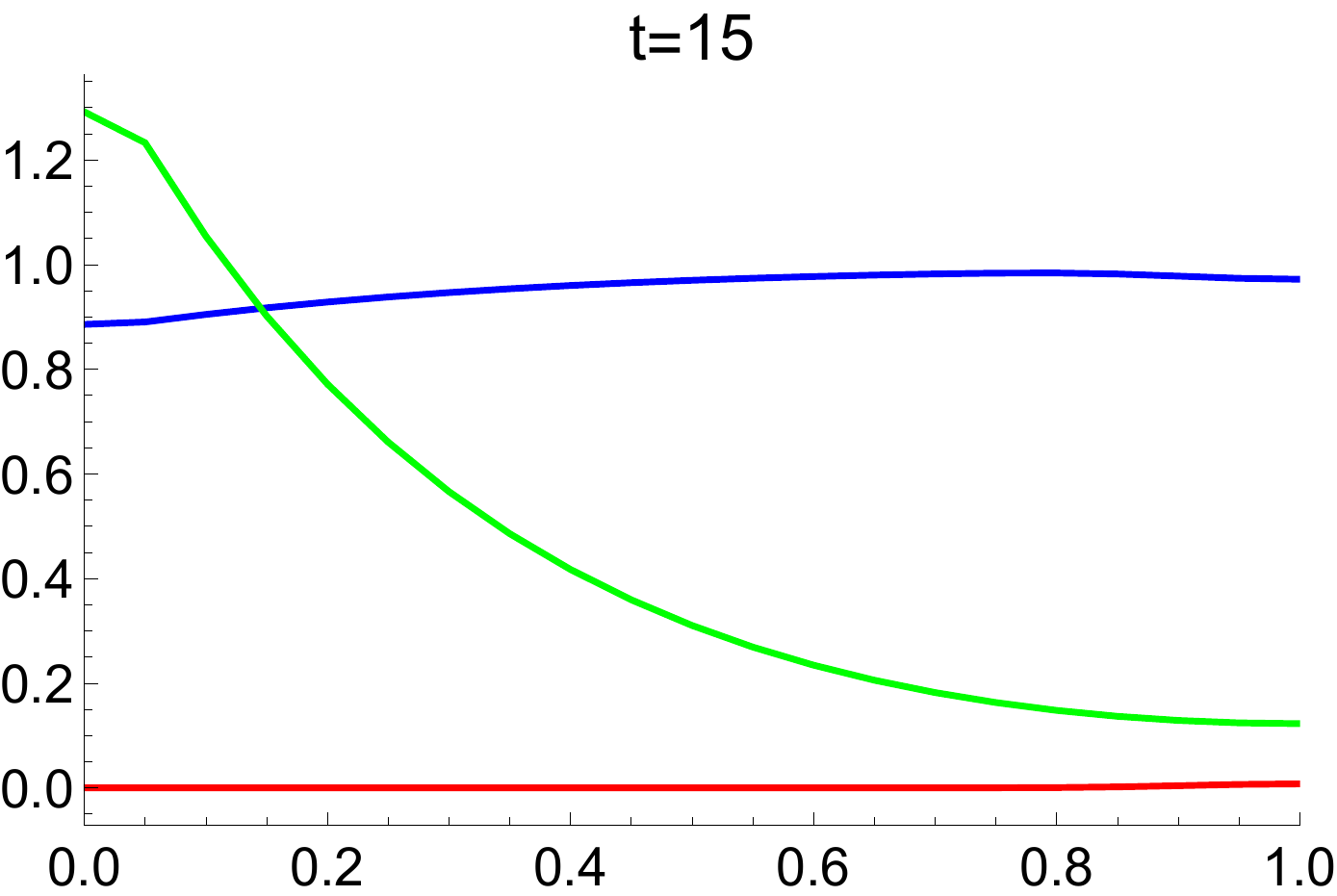}
\caption{Results of Simulation 3, with $\Omega = [0,L]= [0,1]$. Plots of model solutions $A(x,t)$ (cancer cells, red), $N(x,t)$ (normal cells, blue) and $D(x,t)$ (chemotherapeutic drug concentration, green) at time points $t=0,3,6,9,12,15$. See Table \ref{tableSims} for parameter values used here. At time $t=0$, the tumor is spread trough the tissue, and as chemotherapy is applied ($t>0$), the tumor cells are reduced and in the entire tissue. In comparison with Simulation 2, the tumor extinction is reached because the drug infusion rate $\mu$ is increased here. Within the entire tissue, the removal of tumor cells release the normal tissue to recover and grow.}
\label{fig3}
\end{figure}

\begin{figure}[!htb]
\centering
\includegraphics[width=0.32\linewidth]{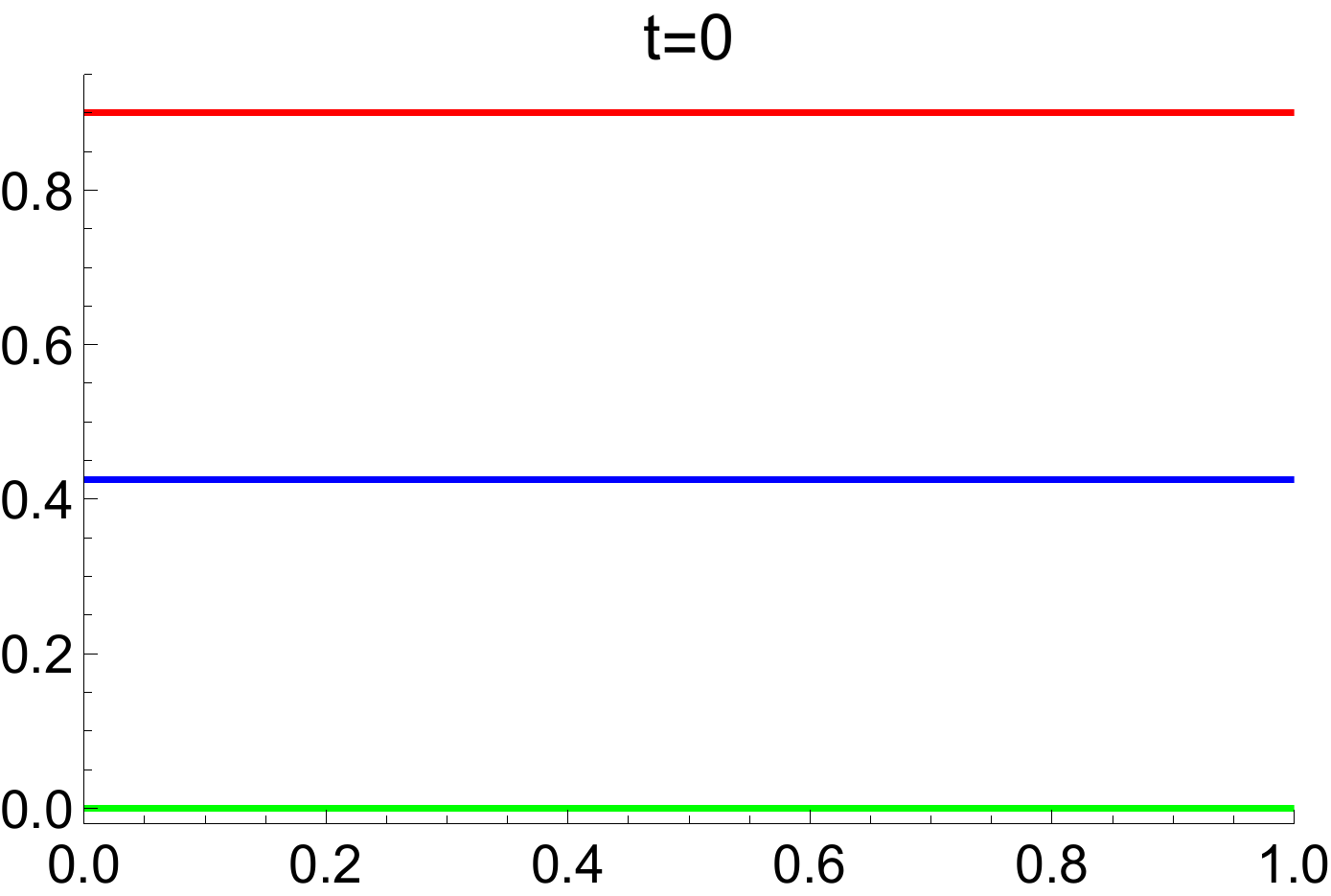}
\includegraphics[width=0.32\linewidth]{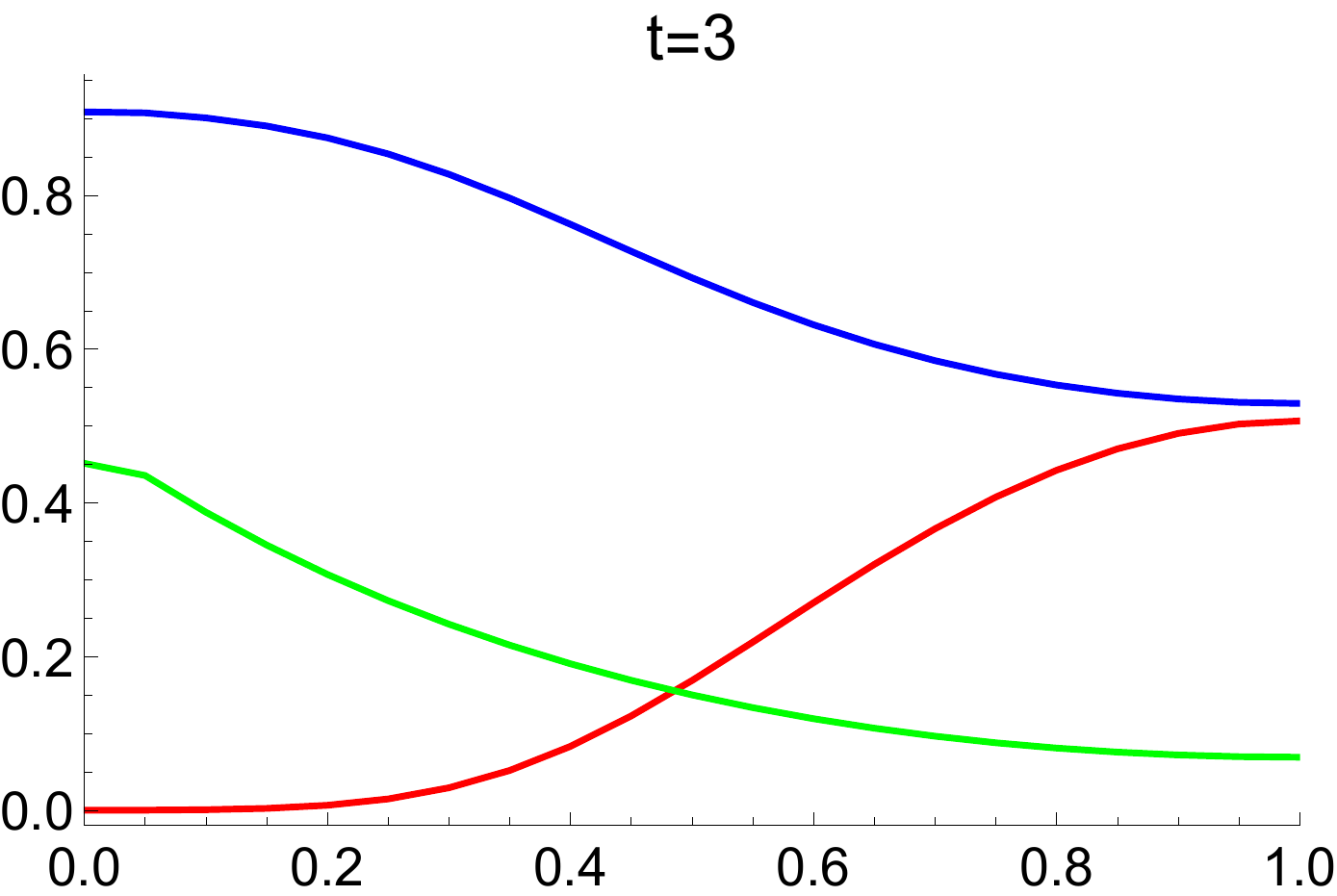}
\includegraphics[width=0.32\linewidth]{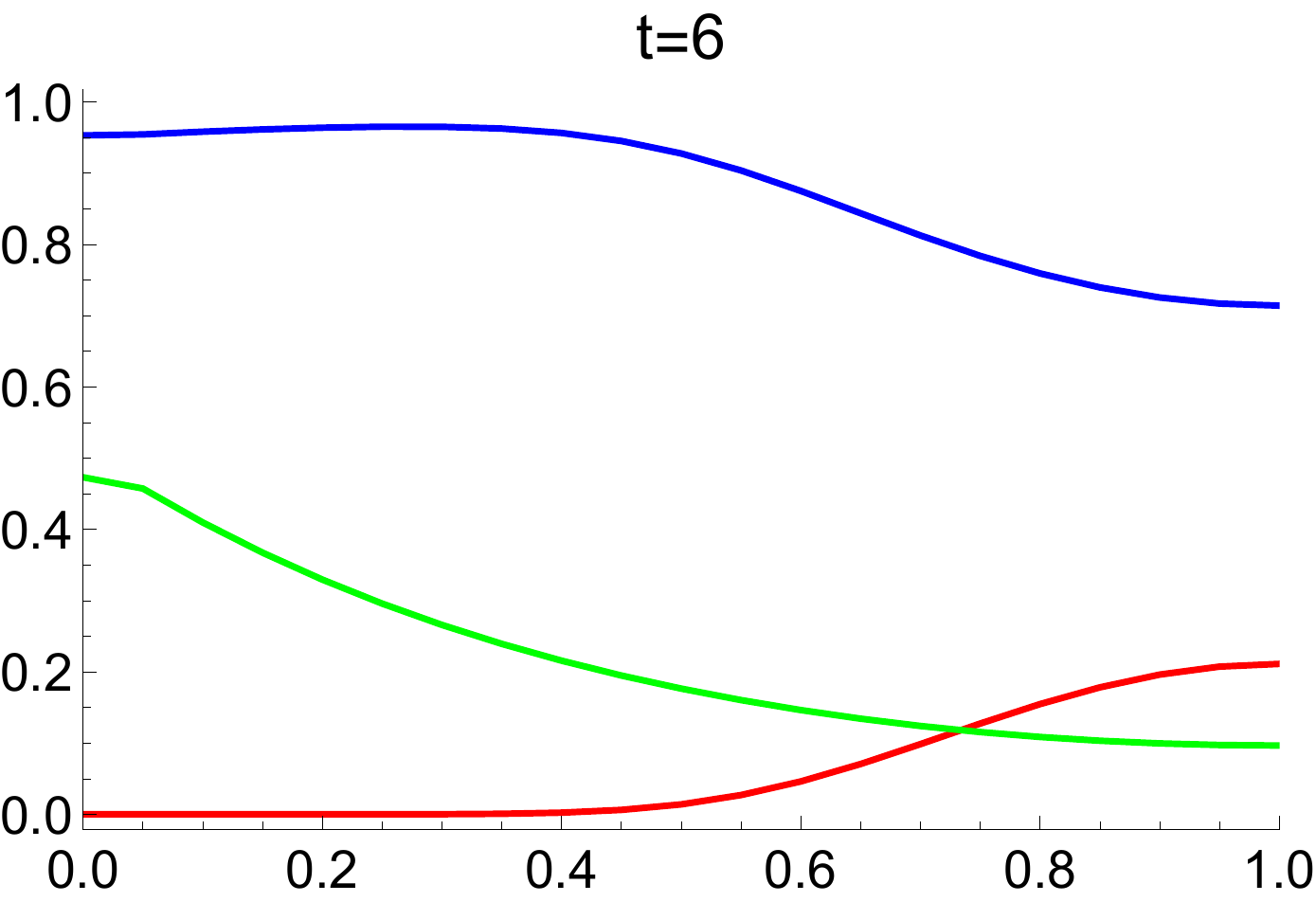}
\includegraphics[width=0.32\linewidth]{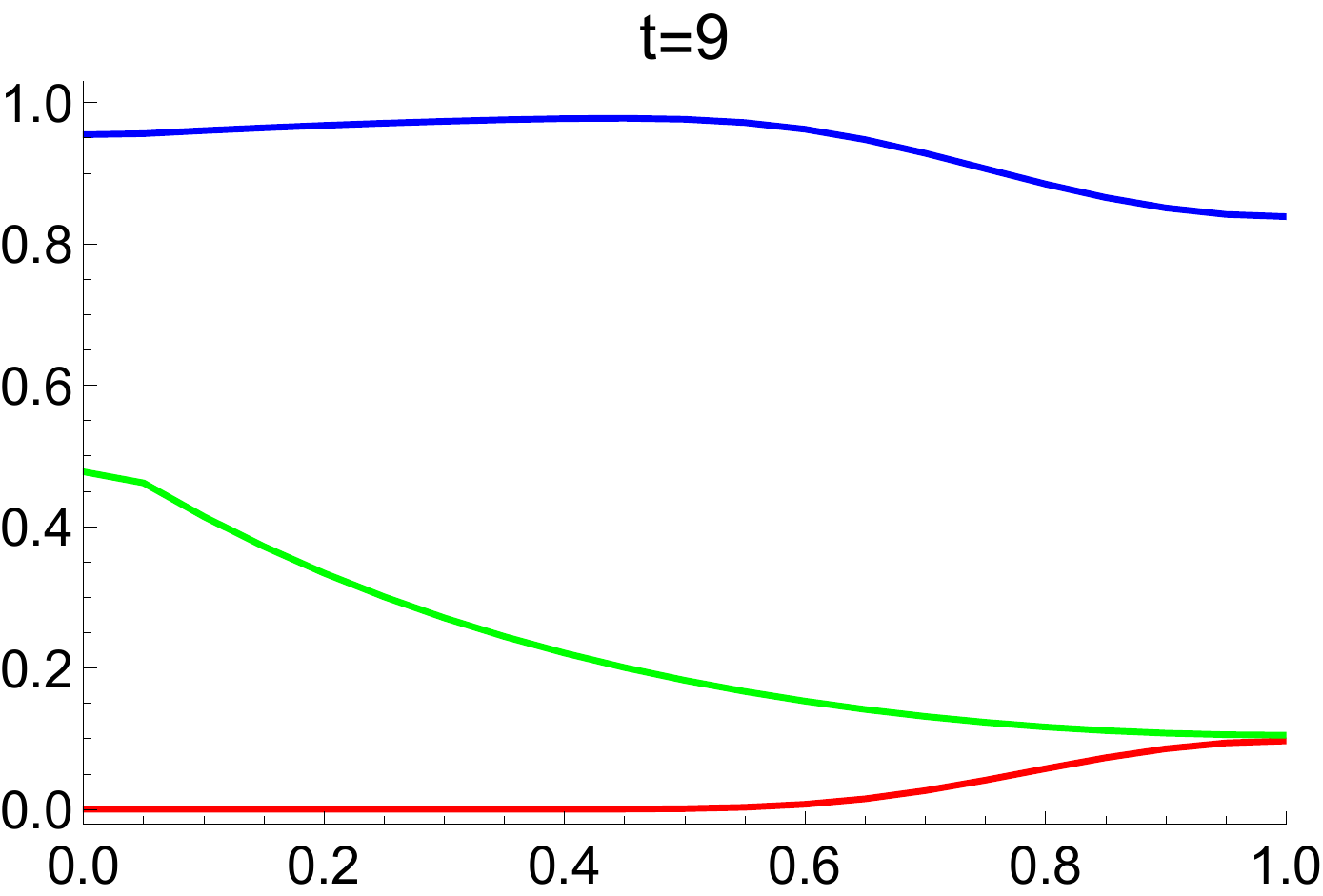}
\includegraphics[width=0.32\linewidth]{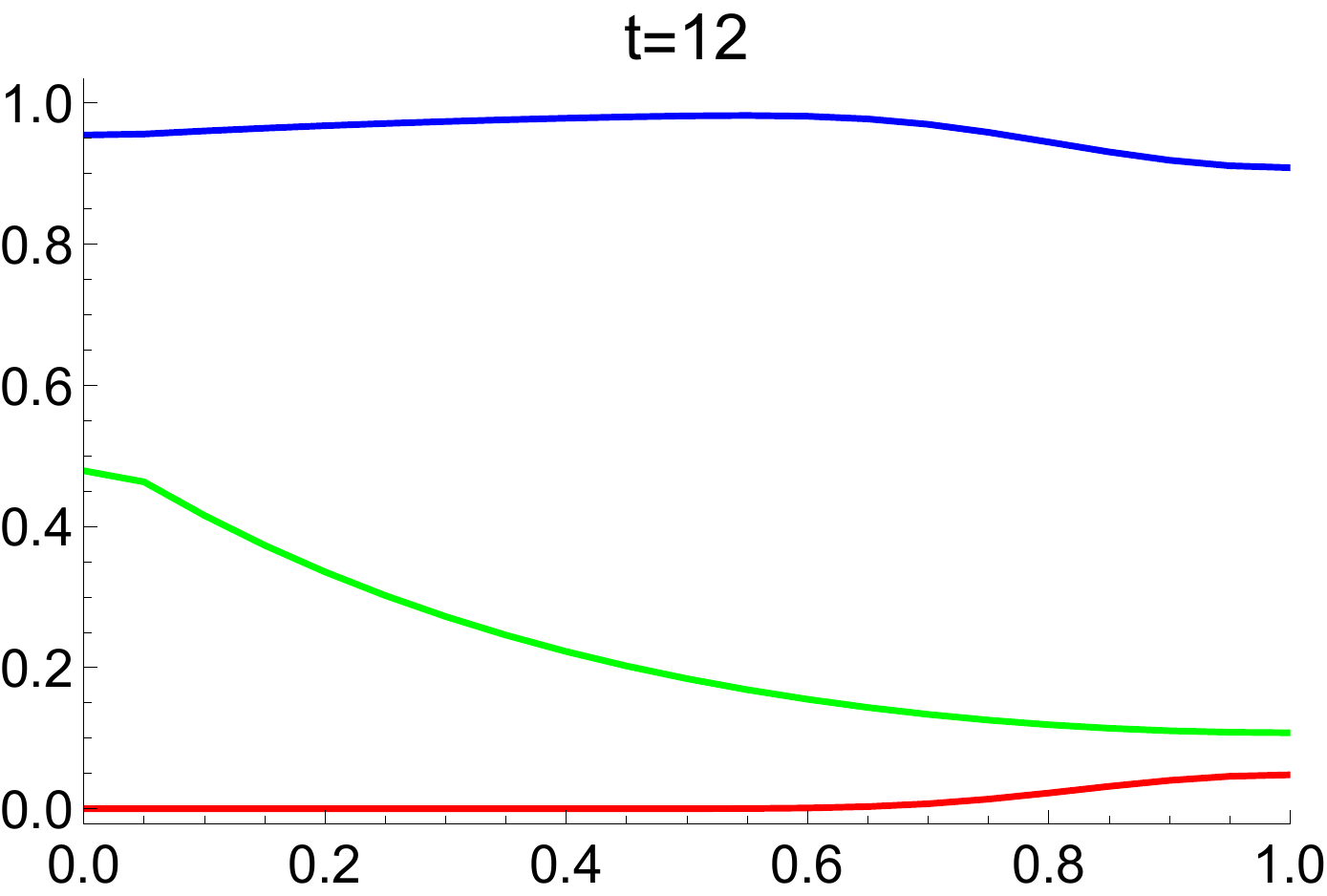}
\includegraphics[width=0.32\linewidth]{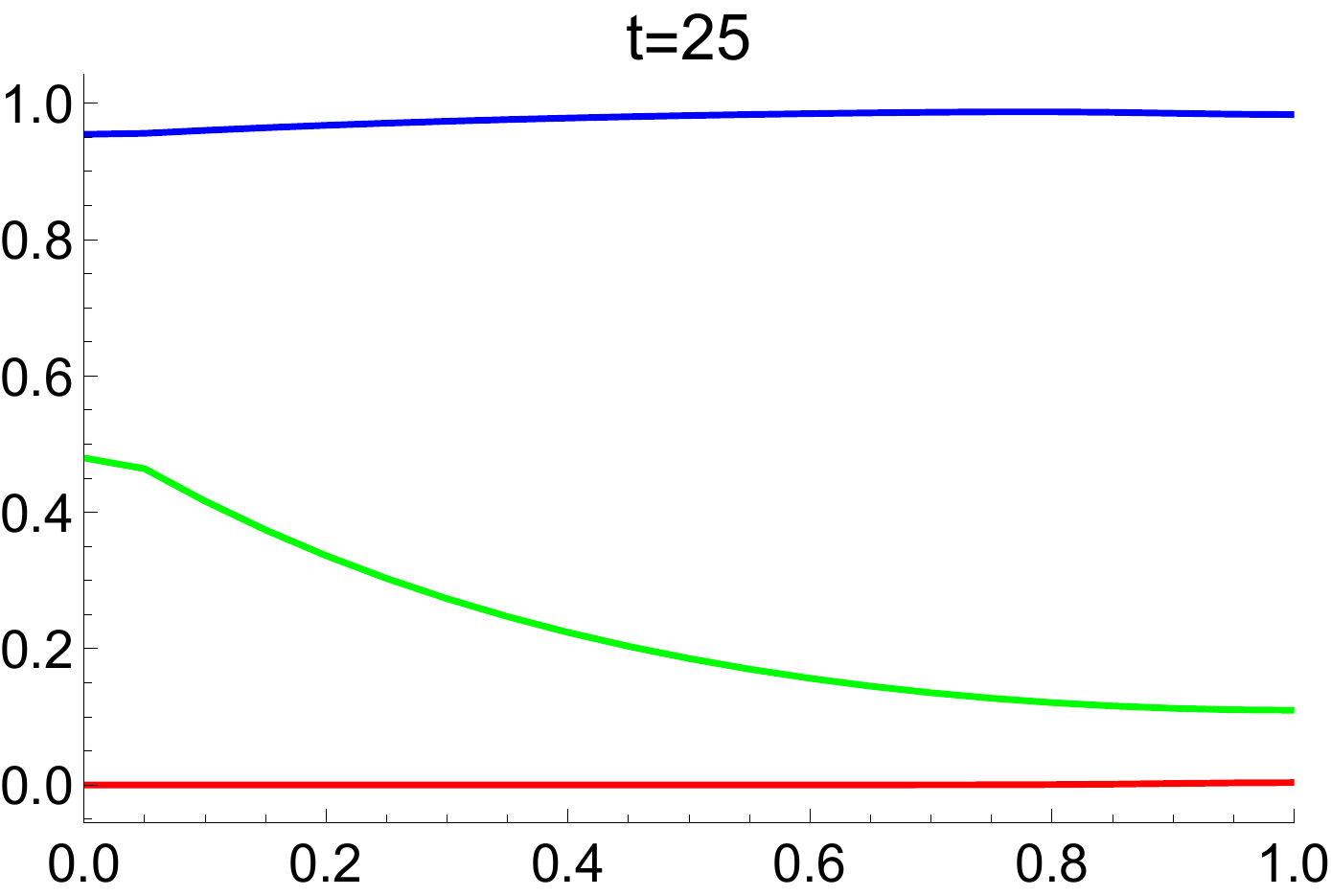}
\caption{Results of Simulation 4, with $\Omega = [0,L]= [0,1]$. Plots of model solutions $A(x,t)$ (cancer cells, red), $N(x,t)$ (normal cells, blue) and $D(x,t)$ (chemotherapeutic drug concentration, green) at time points $t=0,3,6,9,12,25$. See Table \ref{tableSims} for parameter values used here. At time $t=0$, the tumor is spread trough the tissue, and as chemotherapy is applied ($t>0$), the tumor cells are reduced and in the entire tissue. In comparison with Simulation 2, the tumor extinction is reached because the drug diffusion coefficient $\sigma$ is increased here.}
\label{fig4}
\end{figure}

\begin{figure}[!htb]
\centering
\includegraphics[width=0.32\linewidth]{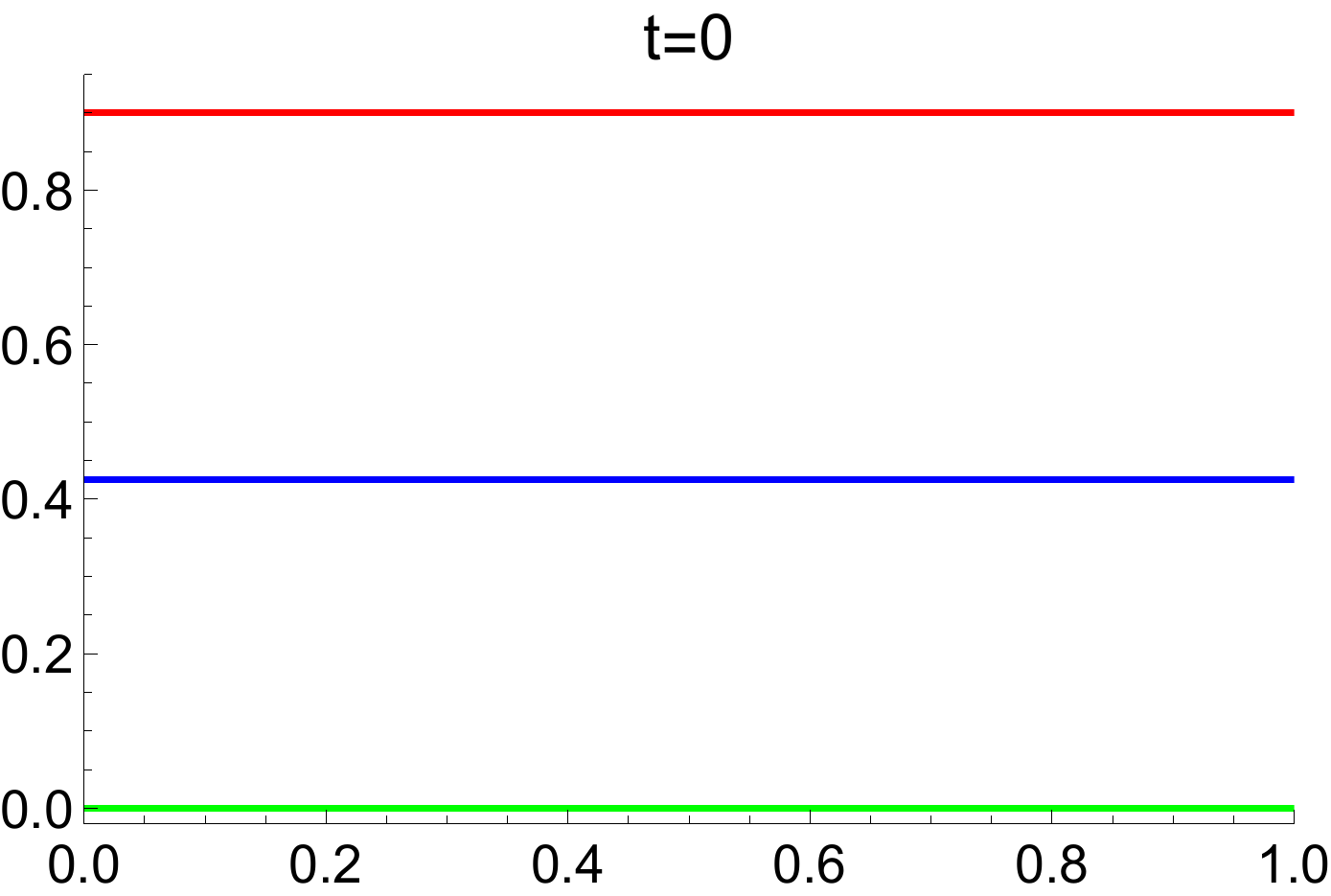}
\includegraphics[width=0.32\linewidth]{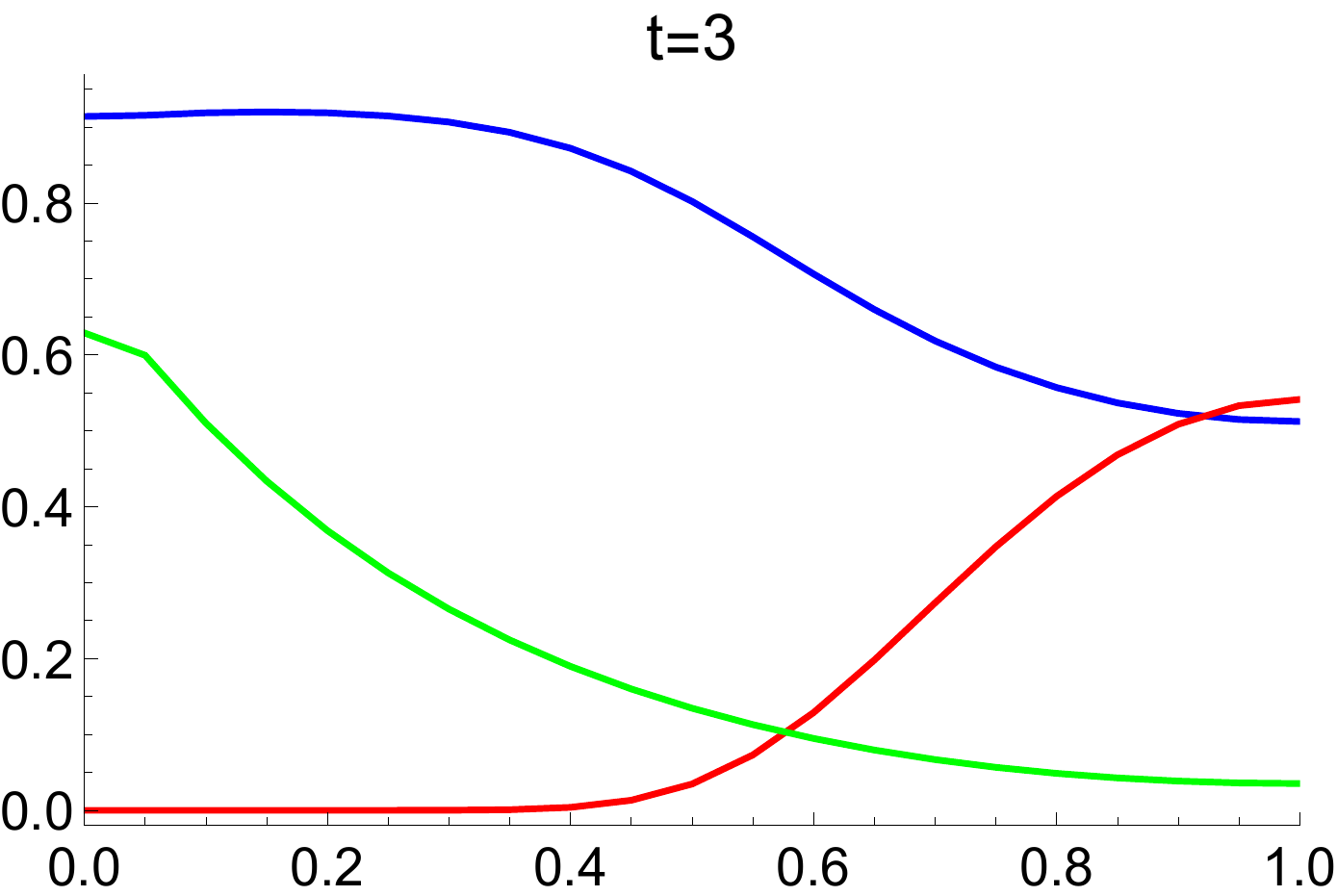}
\includegraphics[width=0.32\linewidth]{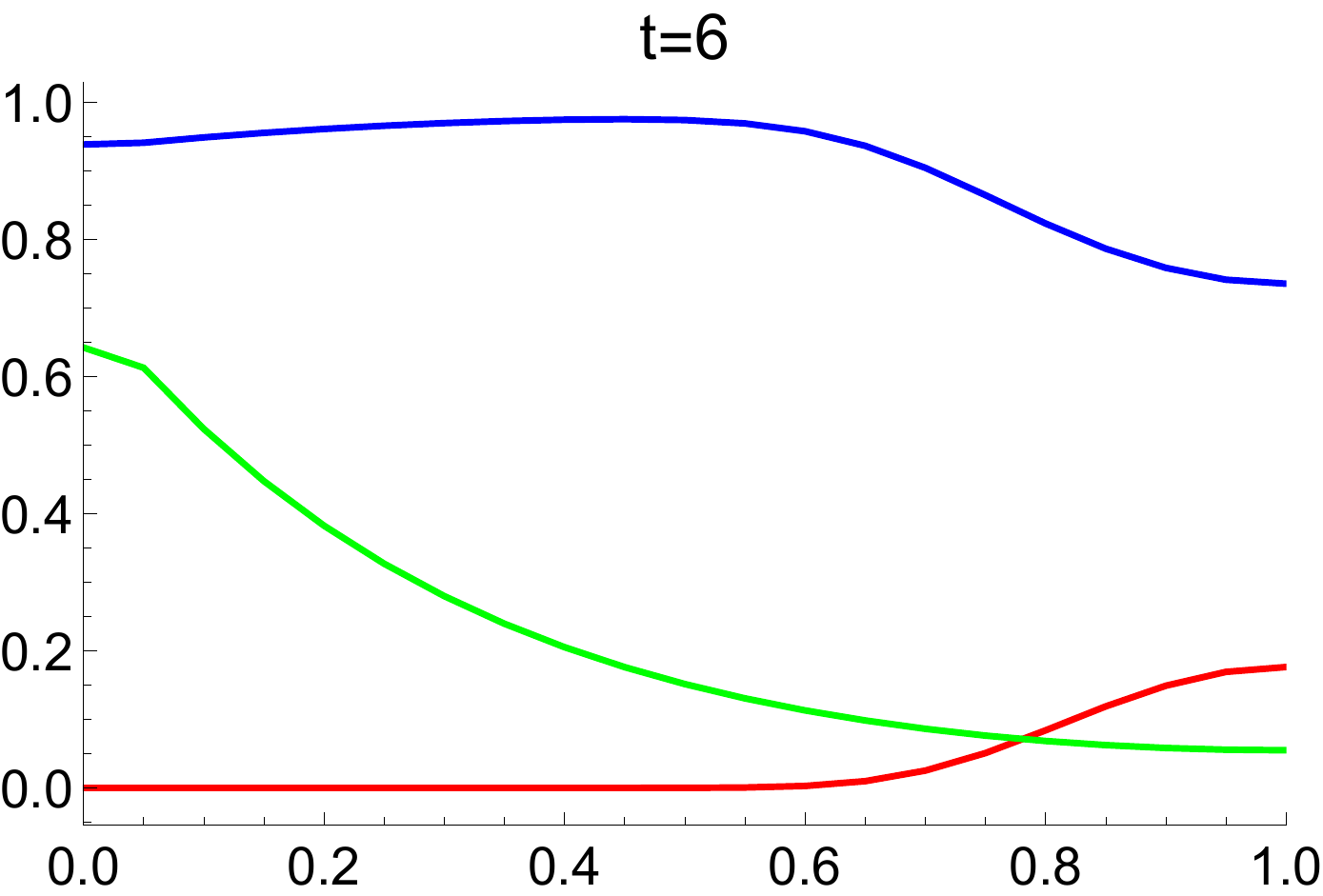}
\includegraphics[width=0.32\linewidth]{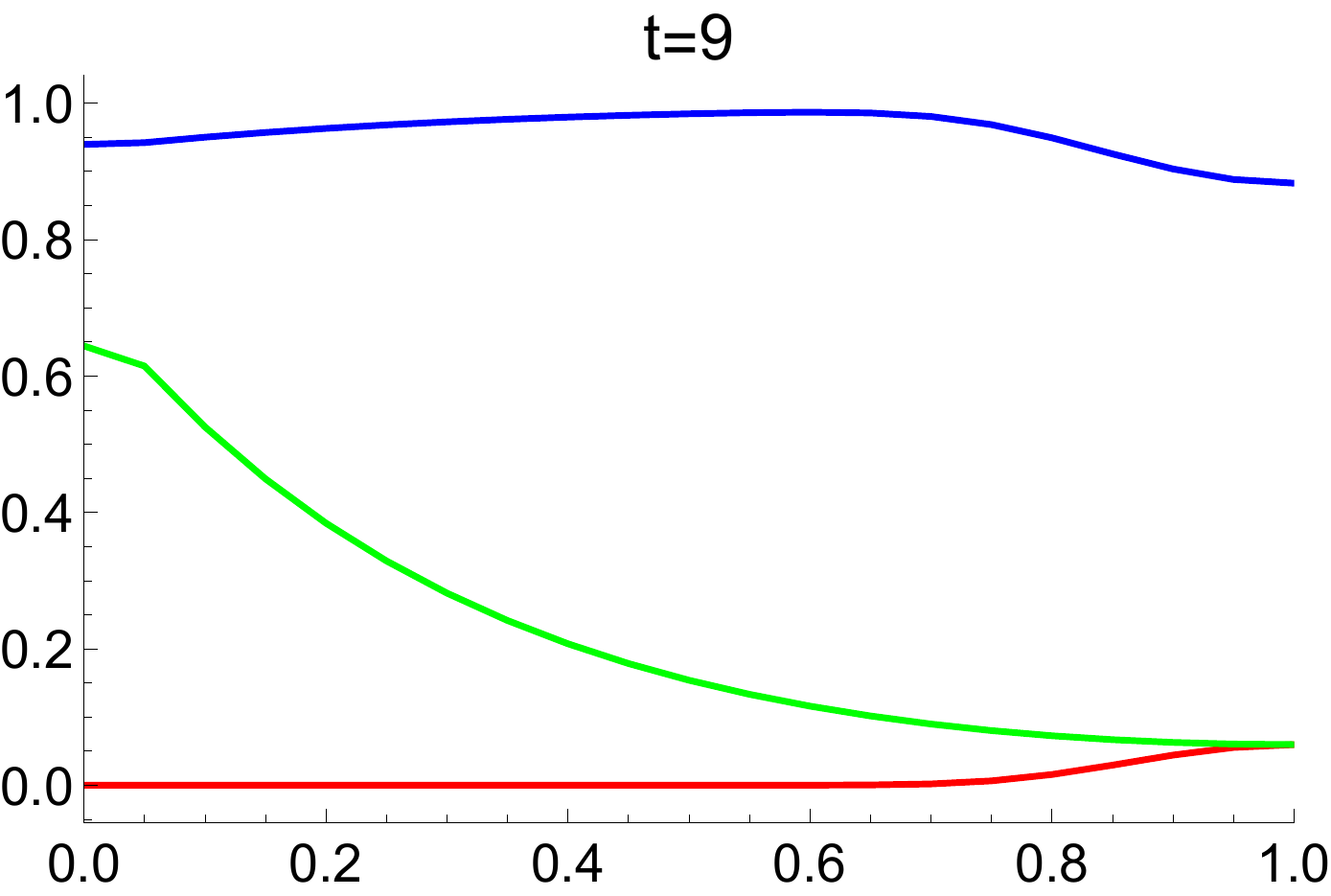}
\includegraphics[width=0.32\linewidth]{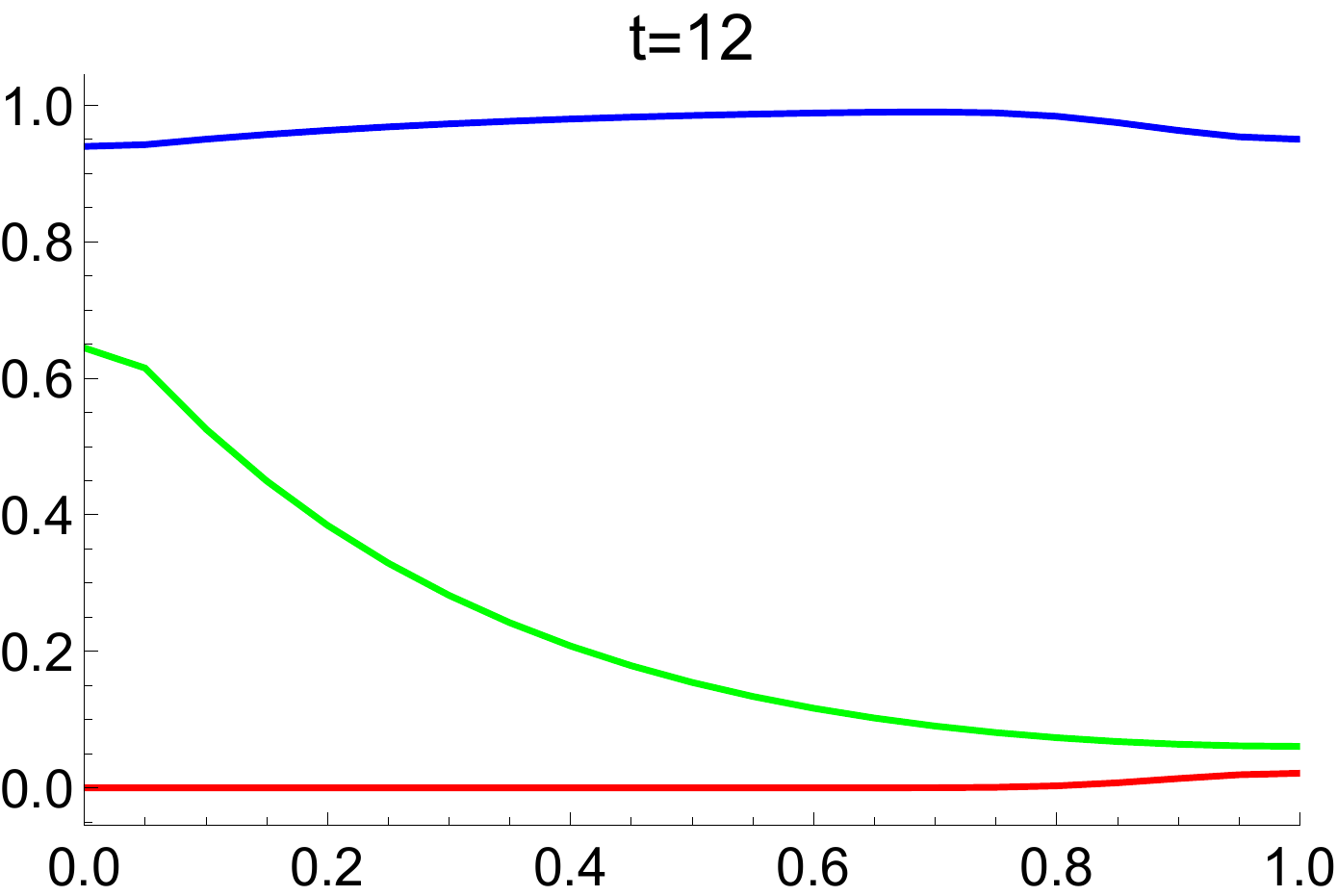}
\includegraphics[width=0.32\linewidth]{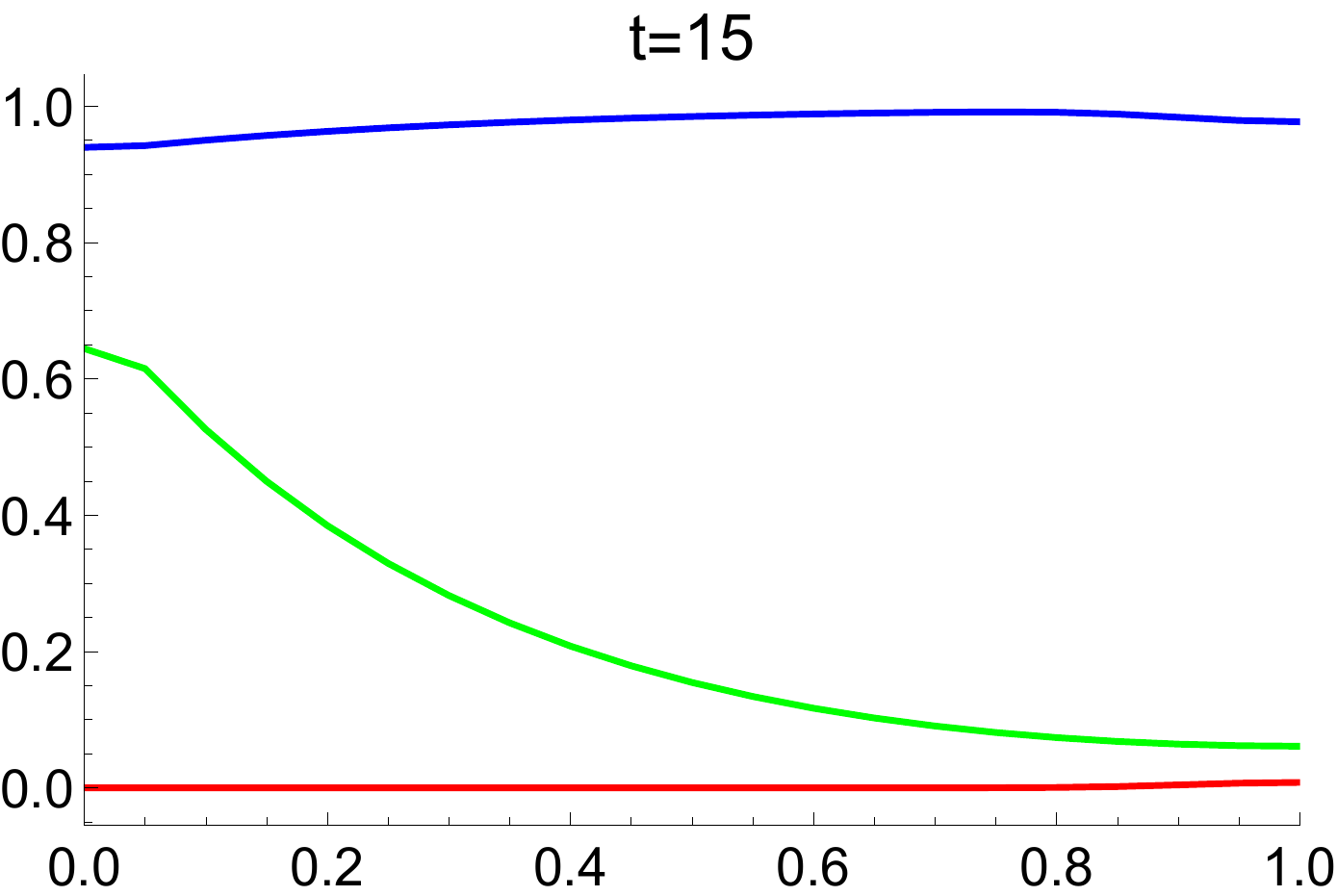}
\caption{Results of Simulation 5, with $\Omega = [0,L]= [0,1]$. Plots of model solutions $A(x,t)$ (cancer cells, red), $N(x,t)$ (normal cells, blue) and $D(x,t)$ (chemotherapeutic drug concentration, green) at time points $t=0,3,6,9,12,15$. See Table \ref{tableSims} for parameter values used here. At time $t=0$, the tumor is spread trough the tissue, and as chemotherapy is applied ($t>0$), the tumor cells are reduced and in the entire tissue. In comparison with Simulation 2, the tumor extinction is reached because the chemotherapy toxicity against tumor cells, $\alpha_A$, was increased.}
\label{fig5}
\end{figure}

In order to make easier to illustrate the model dynamics, we present the results of next simulations in a one-dimensional domain $\Omega = [0, 1]$. In Simulation 2, we use the same parameters values used in Simulation 1 (see Table \ref{tableSims}), but increase the chemotherapy toxicity $\alpha_A$ and move the blood vessel to the left side of the tissue, $\omega=[0,0.1]$. Although the tumor cells in the vicinity of the blood vessel are extinct, the chemotherapy is still not able to eliminate the distant tumor cells (Figure \ref{fig2}). Thus, we observe tumor persistence in the long-term. In Simulation 3, we keep the parameters as in Simulation 2, but increase the chemotherapy infusion rate $\mu$ (mimicking a higher dose). We observe that the tumor cells are extinct in the entire tissue (Figure \ref{fig3}). In Simulation 4, we illustrate other mechanism to achieve tumor extinction: instead of increasing drug dose, we adopt the parameter values of Simulation 2, but increase the drug diffusion $\sigma$, so that it is capable to spread over the entire tissue and effectively eliminate all tumor cells  (Figure \ref{fig4}). Finally, in Simulation 5, we also adopt the parameter values of Simulation 2, but increase the chemotherapy toxicity against tumor cells $\alpha_A$. This also leads to tumor extinction (Figure \ref{fig5}). An advantage of the strategies adopted in Simulations 4 and 5, in comparison with Simulation 3 (increasing dose), is that the former lead to less side effects. Simulation 3 describes the use of a drug which spreads faster, while Simulation 5 illustrates the use of a more potent and specific drug, which targets more tumor cells but not more normal cells ($\alpha_N$ was not changed). Taken together, these simulations and the different outcomes observed for different parameter values confirm the ability of the model to consistently describe tumor chemotherapy and illustrate the potential of mathematical models to provide testable hypothesis that could be studied together with clinicians in order to achieve better results in the treatment of cancer.

\bibliographystyle{amsplain}

\end{document}